\documentclass[final,leqno,showlabel]{siamltex}
\usepackage{amsmath}
\usepackage{amssymb}
\usepackage{cases}
\usepackage{enumerate}
\usepackage{graphicx}
\usepackage{indentfirst}
\usepackage{doi}
\usepackage{cleveref}
\usepackage[usenames]{color}
\usepackage{wrapfig}
\usepackage{stmaryrd}
\usepackage{array} 
\usepackage[caption=false]{subfig} 
\SetSymbolFont{stmry}{bold}{U}{stmry}{m}{n}

\usepackage[notref,notcite]{showkeys}
\usepackage{mathrsfs}
\usepackage{pgfplots}
\usepackage{ulem}
\usepackage{caption}
\usepackage{enumitem}
\usepackage{verbatim}
\pgfplotsset{compat=1.7}

\newtheorem{thm}{Theorem}[section]
\newtheorem{lem}[thm]{Lemma}
\newtheorem{rem}[thm]{Remark}

\renewcommand{\theequation}{\arabic{section}.\arabic{equation}}

\newcommand*{\vv}[1]{\vec{\mkern0mu#1}}

\newcommand{\norm}[1]{\Vert#1\Vert}

\newcommand{\bR}{{\mathbb R}}
\newcommand{\mT}{\mathscr{T}}
\newcommand{\mQ}{\mathcal{Q}}

\newcommand{\bkap}{{\overline{\varkappa}}}
\newcommand{\Gt}{{\Gamma(t)}}
\newcommand{\Gm}{{\Gamma^m}}
\newcommand{\bV}{\mathbb{V}}
\newcommand{\bK}{\mathbb{K}}

\newcommand{\dH}{{\rm d}\mathcal{H}}

\newcommand{\mH}{\mathcal{H}}

\newcommand{\id}{{\rm id}}
\newcommand{\dd}[1]{\frac{\rm d}{{\rm d}#1}}
\newcommand{\ddt}{\dd{t}}
\newcommand{\nn}{\nonumber}
\newcommand{\ttau}{\Delta t}
\newcommand{\ipd}[1]{\bigl(#1\bigr)}
\newcommand{\nabs}{\nabla_{\!s}}

\newcommand{\mat}[1]{\uuline{#1}\rule{0pt}{0pt}}

\title{An energy-stable parametric finite element  method for Willmore flow with normal-tangential velocity splitting  }

\author{Harald Garcke\thanks{Fakult{\"a}t f{\"u}r Mathematik, Universit{\"a}t Regensburg, 
93040 Regensburg, Germany (harald.garcke@ur.de)}
\and Robert N\"urnberg\thanks{Dipartimento di Mathematica, Universit\`a di Trento,
38123 Trento, Italy (robert.nurnberg@unitn.it)}
\and Quan Zhao\thanks{School of Mathematical Sciences, University of Science and Technology of China, 230026 Hefei, Anhui, China (quanzhao@ustc.edu.cn)}
}

\date{}
\begin{document}

\maketitle

\begin{abstract}
We propose and analyze an energy-stable fully discrete parametric approximation for Willmore flow of hypersurfaces in two and three space dimensions. We allow for the presence of spontaneous curvature effects and for open surfaces with boundary. The presented scheme is based on a new geometric partial differential equation  (PDE) that combines an evolution equation for the mean curvature with a separate equation that prescribes the tangential velocity.  The mean curvature is used to determine the normal velocity within the gradient flow structure, thus guaranteeing an unconditional energy stability for the discrete solution upon suitable discretization.  We introduce a novel weak formulation for this geometric PDE, in which different types of  boundary conditions can be naturally enforced. We further discretize the weak formulation to obtain a fully discrete parametric finite element method, for which well-posedness can be rigorously shown.  Moreover, the constructed scheme admits an unconditional stability estimate in terms of the discrete energy. Extensive numerical experiments are reported to showcase the accuracy and robustness of the proposed method for computing Willmore flow of both curves in $\bR^2$ and surfaces in $\bR^3$.
%
\end{abstract}




\pagestyle{myheadings} \markboth{H. Garcke, R. N\"urnberg and Q.~Zhao}
{An energy-stable parametric FEM for Willmore flow }

\section{Introduction} \label{sec:intro}

Geometric functionals involving curvature quantities play an important role in differential geometry, applied mathematics and a variety of physical models. The simplest example is the Willmore energy,  which is defined as the surface integral of the square of the mean curvature  \cite{Willmore93}.  The $L^2$-gradient flow of the Willmore functional gives rise to the so-called Willmore flow,  which is a fourth-order geometric PDE and has a highly nonlinear structure. Willmore flow is a fundamental geometric evolution equation with widespread practical applications, but its intrinsic complexity pose significant challenges for both mathematical analysis and numerical approximations. We refer the reader to \cite{Kuwert01willmore, DziukKS02, KS02, Simonett05willmore, Blatt09singular, Schlierf25spont} and the references therein for relevant analytical results. In this article, we will mainly focus on the numerical approximation of Willmore flow. 

One of the main challenges in the design of numerical methods for Willmore flow is the ability to mimic the gradient flow structure on the (fully) discrete level.
Different approaches have been employed for Willmore flow of hypersurfaces in $\bR^3$ in the literature, including the level set method \cite{Droske04level}, the phase-field method \cite{Du04phase,FrankenRW13,Bretin15phase,RumpfSS25preprint} and front-tracking methods \cite{Mayer02numerical, ClarenzDDRR04, Rusu05, BGN08willmore,Dziuk08,pwfade,pwfopen,OlischlagerR09,BalzaniR12,KovacsLL21,Duan2021high}. However, in general it does not seem possible to derive energy stability estimates for these methods. In what follows we briefly discuss what is known regarding energy stability for these methods. The approaches in \cite{FrankenRW13,RumpfSS25preprint,OlischlagerR09,BalzaniR12} are based on minimizing movement schemes with nested variational time discretizations. As such, they naturally inherit energy stability. With regards to numerical methods that discretize a PDE formulation, Dziuk \cite{Dziuk08} introduced a first parametric finite element method (FEM), for which an energy stability bound can be established in the continuous-in-time semidiscrete setting. This approach was further extended by Barrett, Garcke and N\"urnberg (BGN) in \cite{pwfade, pwfopen} by allowing for a tangential distribution of points, leading to a parametric FEM that not only satisfies the energy decay property on the semidiscrete level, but also preserves the mesh quality of the evolving discrete surfaces. To the best of our knowledge, our manuscript presents the first fully discrete parametric FEM for the Willmore flow of surfaces that is energy stable.

As far as parametric finite element methods are concerned, a second crucial challenge is maintaining the mesh quality of the evolving surface in order to accurately approximate the geometric flow.  For example, discretizing a purely normal flow \cite{DziukKS02, ClarenzDDRR04, Dziuk08, Rusu05, KovacsLL21, Duan2021high} will in general lead to poor surface elements, where mesh smoothing is often necessary to improve the mesh quality.  Of course, such heuristics invalidate any theoretical results shown for the original method, and may pollute the numerical solutions.   A natural remedy  is to incorporate an appropriate tangential velocity into the geometric evolution equations \cite{BGN07, DeTurck17, Duan24new, Hu22evolving}.  We refer the reader to \cite{Barrett20} for the `BGN' approach, which has already been applied to many curvature-driven interface evolution problems. The BGN method relies on a curvature identity (see \eqref{eq:curidenty} below), that under semidiscretization guarantees the equidistribution of mesh points for curve evolutions in $\bR^2$, and so-called conformal polyhedral surfaces in $\bR^3$. The method presented in this paper splits the treatment of normal and tangential motion, and so the desired tangential velocity can be chosen freely. 

Let us mention that our work is greatly motivated by a recent breakthrough due to Bao and Li \cite{BaoL25}, who were able to introduce an energy stable fully discrete parametric FEM for the Willmore flow of closed planar curves. In particular, they considered a system of PDEs that includes an evolution equation for curvature. Their method also incorporates an implicitly defined tangential motion, which is not always benign. Hence mesh redistributions are still required for complex curve evolutions. Recently, in \cite{GNZ25willmore}, the present authors adapted the approach from \cite{BaoL25} to include the BGN tangential motion to yield a numerical method for the Willmore flow of closed planar curves, that is unconditionally energy stable and satisfies an equidistribution property.
Recall that the novel idea of incorporating evolution equations for geometric quantities, like curvature, into the system of PDEs was previously exploited in \cite{KovacsLL21} in the context of a convergence proof for a finite element method for Willmore flow of surfaces. We note that, very recently, a generalization of the work \cite{BaoL25} to the Willmore flow of surfaces has been presented in \cite{BaoLW2025}.

Motived by our recent work for closed curves in the plane \cite{GNZ25willmore}, it is the main aim of this work to propose an energy-stable approximation for Willmore flow of hypersurfaces in a unified framework. Additionally, we will include spontaneous curvature effects and the case of open hypersurfaces with boundary. One of the key ideas in our approach is the introduction of a new geometric PDE for Willmore flow, where an evolution equation of the mean curvature and a curvature identity are included. We then propose a novel weak formulation and discretize the curvature evolution equation in a way that is inspired by the arbitrary Lagrangian--Eulerian method \cite{GNZ24ALE}. This results in a linear method which not only has good mesh property but also satisfies an unconditional energy stability bound. The main contribution of this work can be summarized as follows: 
\begin{itemize}
	\item To the best of our knowledge, this is the first parametric FEM for the Willmore flow of surfaces which satisfies an energy decay property for the fully discrete solutions. 
    \item The method leads to a linear system of equations at each time level.
	\item The approximation of the gradient flow structure is split from applying tangential motion, which allows for flexibility in choosing a tangential velocity without affecting the stability estimate. 
	\item The method can naturally include different types of boundary conditions arising for surfaces with boundary. 
\end{itemize}
For simplicity, in our presentation we will mainly focus on the BGN tangential velocity \cite{Barrett20}. Moreover, for surfaces with boundary we concentrate on two types of boundary conditions: clamped and Navier boundary conditions.

The rest of the article is organized as follows. We begin in Section~\ref{sec:mathform} with the mathematical formulation of Willmore flow. There we introduce a new geometric PDE for describing the considered flow with different boundary conditions. Next in Section~\ref{sec:weak}, we propose a novel weak formulation for this geometric PDE and prove the stability estimate within the weak formulation. In Section~\ref{sec:pfem}, we present our linear parametric finite element approximation and prove its well-posedness and unconditional energy stability. We further report several numerical examples to confirm the convergence, stability and good mesh properties of our introduced scheme in Section~\ref{sec:num}. 

\section{Mathematical formulations}\label{sec:mathform}

Let $(\Gt)_{t\in[0,T]}$ be an evolving 
hypersurface in $\bR^d$  with its parameterization given by
\begin{equation}\label{eq:para}
\vec x(\cdot, t): \Upsilon\times[0,T]\mapsto\bR^d,\qquad d\in\{2,3\},
\end{equation}
where $\Upsilon\subset\bR^d$ is a fixed oriented 
reference manifold, with or without boundary. We set 
\begin{equation*}
	\mathcal{G}_T = \cup_{t\in[0,T]}(\Gamma(t)\times\{t\}).
\end{equation*}
The induced velocity of $\Gt$ by this parameterization is defined on $\mathcal{G}_T$ as
\begin{equation*}
\mathscr{\vv V}(\vec x(\vec\rho,t), t) = \partial_t\vec x(\vec\rho,t)\qquad\forall(\vec\rho,t)\in\Upsilon\times [0,T],
\end{equation*}
while its normal velocity is defined as 
$\mathscr{V} = \mathscr{\vv V} \cdot \vec\nu$,
where $\vec\nu$ is the unit normal to $\Gt$.
The mean curvature of $\Gt$ is given by 
\begin{equation*}
	\varkappa = -\nabs\cdot\vec\nu\qquad\mbox{on}\quad\Gt,
\end{equation*}
where $\nabs$ is the surface gradient operator. Our sign convention is such that $\varkappa=-(d-1)$ for the unit sphere with outer normal. Let $\Delta_s=\nabs\cdot\nabs$ be the Laplace-Beltrami operator. We then have the curvature identity 
\begin{equation}
	\varkappa\,\vec\nu = \Delta_s\vec\id\qquad\mbox{on}\quad \Gt,\label{eq:curidenty}
\end{equation}
where  $\vec\id$ is the identity function in $\bR^d$.

Given $f\in C^1(\mathcal{G}_T)$, we define the material derivative
\begin{equation*}
	\partial_t^\circ f = \ddt f(\vec x(\vec\rho,t), t)\quad\forall(\vec\rho,t)\in\Upsilon\times [0,T],
\end{equation*}
 that follows the parameterization \eqref{eq:para}. The normal time derivative of $f$ is defined as
\begin{equation} \label{eq:normalf}
	\partial_t^\square f = \partial_t^\circ f-\mathscr{\vv V}\cdot\nabs f\qquad\mbox{on}\quad\Gt,
\end{equation}
which measures the change of $f$ on the moving surface in the normal direction. The following lemma provides the normal time derivative of the unit normal and the mean curvature. For the proof we refer to \cite[Lemmas~37,38,39]{Barrett20}.
\begin{lem} It holds that
\begin{subequations}
\begin{alignat}{2}\label{eq:timenv}
		\partial_t^\square\vec\nu & = -\nabs\mathscr{V}\quad &\mbox{on}\quad\Gt,\\
		\partial_t^\square\varkappa &= \Delta_s\mathscr{V}+\mathscr{V}|\nabs\vec\nu|^2\quad &\mbox{on}\quad\Gt,\label{eq:timcurva}
\end{alignat}
\end{subequations}
where 
$\nabs\vec\nu$ is the Weingarten map,  and $|\mat{A}|^2={\rm tr}(\mat{A}\,\mat{A}^T)$ is the Frobenius norm for any matrix $\mat{A}\in\bR^{d\times d}$. In the case $d=2$, $|\nabs\vec\nu|^2$ reduces to $\varkappa^2$. 
\end{lem}
 
\subsection{The Willmore flow}

\begin{figure}[!htp]
	\centering
	\includegraphics[width=0.45\textwidth]{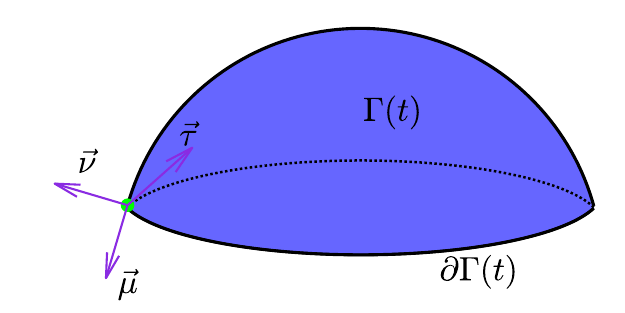}
	\caption{Sketch of $\Gt$ with boundary $\partial\Gt$, as well as the three unit vectors $\vec\tau,\vec\mu$ and $\vec\nu$ which form an orthonormal basis of $\bR^3$.}
	\label{fig:openS}
\end{figure}

For simplicity we define $\vec x(t) = \vec x(\cdot, t)$ and $\varkappa(\cdot, t)=\varkappa(t)$.  In this work, we  consider the energy contribution 
\begin{equation}
E_{\bkap}(\varkappa(t), \vec x(t)) = \frac{1}{2}\int_{\Gamma(t)}(\varkappa-\bkap)^2\,\dH^{d-1},\qquad\bkap\in\bR, \label{eq:Tenergy}
\end{equation}
where the constant $\bkap $ is the so-called spontaneous curvature and  $\dH^{d-1}$ is integration with respect to the $(d-1)$-dimensional Hausdorff measure in $\bR^d$. Here we allow $\Gt$ to be either a closed hypersurface, so that $\partial\Gt=\emptyset$, or an open hypersurface with boundary $\partial\Gt$, as shown in Fig.~\ref{fig:openS}.  Using the results in \eqref{eq:timcurva}, \eqref{eq:transport} as well as integration by parts formula \eqref{eq:interpart},  it is not difficult to obtain that (see \cite[(222)]{Barrett20})
\begin{align}
	&\ddt E_{\bkap}(\varkappa(t),\vec x(t)) = \int_{\Gt}\left[\Delta_s\varkappa + (\varkappa-\bkap)|\nabs\vec\nu|^2 - \frac{1}{2}(\varkappa-\bkap)^2\right]\,\mathscr{V}\dH^{d-1}\nn\\
	&\qquad +\int_{\partial\Gt}\left[\frac{1}{2}(\varkappa-\bkap)^2\,\mathscr{\vv V}\cdot\vec\mu - (\vec\mu\cdot\nabs\varkappa)\,\mathscr{ V} + (\varkappa-\bkap)\,(\vec\mu\cdot\nabs\mathscr{V})\right]\dH^{d-2}, \label{eq:dtE}
\end{align}
where $\vec\mu$ is the outer unit conormal to $\partial\Gt$, and $\dH^{d-2}$ represents integration with respect to the $(d-2)$-dimensional Hausdorff measure in $\bR^d$. 

Taking the $L^2$-gradient flow of  \eqref{eq:Tenergy}  leads to the  Willmore flow
\begin{equation}
\mathscr{V} = -\Delta_s\varkappa -(\varkappa - \bkap)\,|\nabs\vec\nu|^2+\frac{1}{2}(\varkappa-\bkap)^2\varkappa \qquad\mbox{on}\quad\Gt, \label{eq:Willmore}
\end{equation}
together with boundary conditions on $\partial\Gamma(t)$,
in the case of an open hypersurface, that make the
boundary terms in \eqref{eq:dtE} vanish. In this paper, we consider the following two cases
 \begin{itemize}
 	\item [(i)] Navier boundary 
 \begin{subequations}\label{eqn:bd}
 	\begin{equation}
 	\partial\Gt = \partial\Gamma(0),\quad \varkappa =\bkap \quad\mbox{on}\quad\partial_1\Gt;
 		\label{eq:navier}
 	\end{equation}
 	\item  [(ii)] clamped boundary
 	\begin{equation}
 	\partial\Gt = \partial\Gamma(0),\quad \vec\mu\cdot\nabs\mathscr{V}= 0\quad\mbox{on}\quad\partial_2\Gt;
 		\label{eq:clamp}
 	\end{equation}
 \end{subequations}
 \end{itemize}
 where $\partial\Gt =\cup_{i=1}^2\partial_i\Gt$ is a disjoint partitioning of the boundary.  In both cases,  we fix  the boundary, which means that
\begin{equation}
 \mathscr{\vv V}\cdot\vec\mu = 0,\qquad \mathscr{\vv V}\cdot\vec\nu = \mathscr{V}=0\quad\mbox{on}\quad\Gt.\label{eq:fixbd}
 \end{equation}
 This implies that all the boundary terms in \eqref{eq:dtE} vanish for both 
\eqref{eq:navier} and \eqref{eq:clamp}. Hence, \eqref{eq:dtE} shows that the 
geometric flow \eqref{eq:Willmore}, together with the boundary conditions 
\eqref{eqn:bd}, obeys the energy dissipation law
\begin{equation}
	\ddt E_{\bkap}(\varkappa(t), \vec x(t)) = -\int_{\Gt}\mathscr{V}^2\dH^{d-1}\leq 0.
	\label{eq:energylaw}
\end{equation}

We next present a common geometric interpretation of the clamped boundary condition \eqref{eq:clamp}. Let $\vec\tau := \vec\mu \times \vec\nu$ be the unit tangent to $\partial\Gt$, see Fig.~\ref{fig:openS}.
It follows from $\nabs\mathscr{V}\in{\rm span}\{\vec\tau,\vec\mu\}$ and
\eqref{eq:clamp} that $(\nabs\mathscr{V}) \times \vec\tau = \vec 0$.
Hence, on noting that $\vec\tau$ is fixed in time, and on recalling
\eqref{eq:timenv}, we have that
 \begin{equation}
 	\partial_t^\square\vec\mu = \partial_t^\square(\vec\nu\times\vec\tau)=(\partial_t^\square\vec\nu)\times\vec\tau =-(\nabs\mathscr{V})\times\vec\tau=\vv 0\qquad\mbox{on}\quad\partial_2\Gt,\label{eq:clamptau}
 \end{equation}%
which means that the conormal vector $\vec\mu$ does not change in time,
i.e., it is prescribed by the conormal of $\Gamma(0)$.

It is the main aim of this work to devise a fully discrete parametric finite element approximation of the considered  Willmore flow such that a discrete analogue of \eqref{eq:energylaw} holds on the discrete level. 

\subsection{A new geometric PDE}

We next aim to introduce a new geometric PDE for the considered flow. Combining \eqref{eq:normalf} with \eqref{eq:timcurva} yields that
\begin{equation}
\partial_t^\circ\varkappa = \Delta_s\mathscr{V} + \mathscr{V}|\nabs\vec\nu|^2+\mathscr{\vv V}\cdot\nabs\varkappa.\label{eq:curvevo}
\end{equation}
Compare this with \cite[(2.10a)]{KovacsLL21}, where the last term vanishes, since the flow considered there is purely normal, i.e., $\mathscr{\vv V} = \mathscr{V}\,\vec\nu$. In this work we shall consider a flow with an appropriate tangential velocity. The extra term $\mathscr{\vv V}\cdot\nabs\varkappa$ is therefore included to account for  the rate of change of $\varkappa$ induced by the tangent component of $\mathscr{\vv V}$.  

We combine the curvature evolution equation \eqref{eq:curvevo} with the curvature identity \eqref{eq:curidenty}, and introduce the following equivalent system for the Willmore flow on $\mathcal{G}_T$:
\begin{subequations}\label{eqn:refgov}
\begin{align}
\label{eq:rw1}
\mathscr{V} &= -\Delta_s\varkappa - (\varkappa-\bkap)|\nabs\vec\nu|^2 + \frac{1}{2}(\varkappa-\bkap)^2\varkappa,\\
\label{eq:rw2}
 \partial_t^\circ\varkappa &= \Delta_s\mathscr{V} + \mathscr{V}|\nabs\vec\nu|^2+\mathscr{\vv V}\cdot\nabs\varkappa,\\
 \label{eq:rw3}
 \mathscr{\vv V}\cdot\vec\nu & = \mathscr{V},\\
\kappa\,\vec\nu &= \Delta_s\vec\id,
\label{eq:rw4}
\end{align}
\end{subequations}
where $\varkappa$  is the curvature for the gradient flow, while the curvature $\kappa$ acts as a Lagrange multiplier for the side constraint \eqref{eq:rw3}. In the case when $\partial\Gt\neq\emptyset$, we consider the boundary conditions:
	\begin{subequations}\label{eqn:bd2}
	\begin{alignat}{2}
		\label{eq:navier1}
		&\mathscr{\vv V} = \vv0,\quad\mathscr{V} = 0,\quad\varkappa = \bkap \qquad &\mbox{on}\quad\partial_1\Gt,\\
		\label{eq:clamp1}
		&\mathscr{\vv V} = \vv0,\quad\mathscr{V} = 0,\quad\vec\mu\cdot\nabs\mathscr{V} = 0\qquad &\mbox{on}\quad\partial_2\Gt.
	\end{alignat}
	 \end{subequations}
Observe that \eqref{eq:navier1} and \eqref{eq:clamp1} imply the Navier
boundary condition \eqref{eq:navier} and the clamped boundary condition 
\eqref{eq:clamp}, respectively. We note furthermore that we have 
strengthened the conditions \eqref{eq:fixbd} to
$\mathscr{\vv V}=\vv 0$, which means that we do not allow the
parameterization to have a tangential velocity along the fixed boundary line, 
i.e, $\mathscr{\vv V}\cdot\vec\tau=0$.

For the second curvature $\kappa$,  we impose  additionally 
\begin{equation}
	\kappa = \bkap\qquad\mbox{on}\quad\partial\Gt.
\end{equation}
 To complete the system,  we also need the initial curvature  $\varkappa(\cdot,0)$ and the initial surface $\Gamma(0)$.

\begin{rem}
The new system \eqref{eqn:refgov} is the key ingredient that enables an energy-stable approximation and the separation of tangential and normal velocity. The first two equations, \eqref{eq:rw1}–\eqref{eq:rw2}, describe the normal velocity $\mathscr{V}$ and the mean curvature $\varkappa$, reflecting the underlying gradient flow structure. The last two equations, \eqref{eq:rw3}–\eqref{eq:rw4}, can be interpreted as imposing the normal velocity while incorporating the BGN tangential velocity. 
\end{rem}

We recall that the original BGN method \cite{BGN08willmore} for Willmore flow is based on the formulation 
\begin{equation}
 \mathscr{\vv V}\cdot\vec\nu = -\Delta_s\varkappa - (\varkappa-\bkap)|\nabs\vec\nu|^2 + \frac{1}{2}(\varkappa-\bkap)^2\varkappa,\qquad\varkappa\,\vec\nu = \Delta_s\vec\id,\label{eq:BGN1}
\end{equation}
while the stable semidiscrete BGN method from \cite{pwfade} is based on a variational formulation of a system of the form
\begin{equation}
  \mathscr{\vv V} =  \vec f(\Gamma, \vec\varkappa, \vec y), \qquad\vec\varkappa = \Delta_s\vec\id, \qquad \vec y = \vec\varkappa - \bkap\vec\nu.
\label{eq:BGN2}
\end{equation}
For the formulation \eqref{eq:BGN1} it does not appear possible to come up with a variational formulation for the boundary conditions in \eqref{eqn:bd2}, and no stability results are available for its discretizations either. For the formulation \eqref{eq:BGN2}, on the other hand, it is possible to accommodate different boundary conditions in a variational setting, 
see \cite{pwfopen}. However, stability estimates can only be established for the solutions of a continuous-in-time semidiscrete scheme.

\section{Weak formulation}
\label{sec:weak}

In order to propose a weak formulation for the system \eqref{eqn:refgov} with \eqref{eqn:bd2}, we introduce the function spaces on $\Gamma$
\begin{align}
	\bV_{\partial_\alpha}(\Gamma) &:=\bigl\{\chi\in H^1(\Gamma): \chi = 0\quad\mbox{on}\quad\partial_\alpha\Gamma\bigr\},\quad  \alpha=1,2; \nn\\
	\bV_0(\Gamma)&:=\bV_{\partial_1}(\Gamma)\cap\bV_{\partial_2}(\Gamma), \nn 
\end{align}
and the subsets of the space $H^1(\Gamma)$
\begin{align}
\bK_{\partial_\alpha}(\Gamma) &:= \bigl\{\chi\in H^1(\Gamma): \chi = \bkap\quad\mbox{on}\quad\partial_\alpha\Gamma\bigr\},\quad \alpha=1,2;\nn\\
	\bK_{\bkap}(\Gamma) &:=\bK_{\partial_1}(\Gamma)\cap\bK_{\partial_2}(\Gamma). \nn
\end{align}
We further denote by $(\cdot, \cdot)_{\Gamma}$ the  $L^2$-inner product on $\Gamma$.   

We next employ an antisymmetric treatment for the convective term $\mathscr{\vv V}\cdot\nabs\varkappa$ in \eqref{eq:rw2}, which will be beneficial for the stability estimate.
\begin{lem}
Let $\Gamma(t)$ be an evolving hypersurface with $\partial\Gamma(t) =
\partial\Gamma(0)$ for all $t \in [0,T]$. Then it holds that
	\begin{align}
		&\ipd{\mathscr{\vv V}\cdot\nabs\varkappa,~\chi}_\Gt + \frac{1}{2}\bigl(\nabs\cdot\mathscr{\vv V},~[\varkappa-\bkap]\,\chi\bigr)_{\Gt}\nn\\
	&\qquad  =  \frac{1}{2}\mathscr{A}_{\Gt}(\mathscr{\vv V}, \varkappa-\bkap,~\chi)-\frac{1}{2}\bigl(\mathscr{\vv V}\cdot\vec\nu,~[\varkappa-\bkap]\,\varkappa\,\chi\bigr)_{\Gt}\qquad \forall\chi\in H^1(\Gt), \label{eq:anti}
	\end{align}
where $\mathscr{A}_{\Gt}$ is the antisymmetric term defined via
\begin{equation} \label{eq:antisym}
\mathscr{A}_{\Gamma}(\vec\eta, u,v)=\left[\bigl(\vec\eta\cdot\nabs u, ~v\bigr)_{\Gamma}-\bigl(\vec\eta\cdot\nabs v, ~u\bigr)_{\Gamma}\right].
\end{equation}
\end{lem}

\begin{proof}
The following identity holds  
\begin{align}
	(\nabs\varkappa)\,\chi & = \nabs(\varkappa-\bkap)\,\chi\nn\\ 
	&= \frac{1}{2}\nabs(\varkappa-\bkap)\,\chi +\frac{1}{2}\left[\nabs([\varkappa-\bkap]\,\chi) - (\varkappa-\bkap)\nabs\chi\right].
\end{align}
Thus a direct computation yields 
\begin{align}
&\bigl(\mathscr{\vv V}\cdot\nabs\varkappa, ~\chi\bigr)_{\Gt}= \frac{1}{2}\mathscr{A}_{\Gt}(\mathscr{\vv V}, \varkappa-\bkap,~\chi)+\frac{1}{2}\bigl(\mathscr{\vv V},~\nabs[(\varkappa-\bkap)\,\chi]\bigr)_{\Gt}\nn\\
&=\frac{1}{2}\mathscr{A}_{\Gt}(\mathscr{\vv V}, \varkappa-\bkap,~\chi)-\frac{1}{2}\bigl(\mathscr{\vv V}\cdot\vec\nu,~[\varkappa-\bkap]\,\varkappa\,\chi\bigr)_{\Gt}-\frac{1}{2}\bigl(\nabs\cdot\mathscr{\vv V},~[\varkappa-\bkap]\,\chi\bigr)_{\Gt},\nn
\end{align}
where we have used the integration by parts formula \eqref{eq:interpart} and recalled that $\mathscr{\vv V} = 0$ on $\partial\Gt$.  
\end{proof}

We now propose the following weak formulation for the new geometric PDE. Let $\vec x(\cdot,0)\in [H^1(\Upsilon)]^d$   with  $\vec x(\Upsilon, 0)= \Gamma(0)$  and  $\varkappa(\cdot,0)\in \bK_{\bkap}(\Gamma(0))$.  For each $t\in(0,T]$ we find $\vec x(\cdot, t)\in [H^1(\Upsilon)]^d$ for $\Gt=\vec x(\Upsilon,t)$ with $\mathscr{\vv V}(t)\in [\bV_0(\Gt)]^d$, $\mathscr{V}(t)\in \bV_0(\Gt)$, $\varkappa(t)\in\bK_{\partial_1}(\Gt)$ and $\kappa(t)\in \bK_{\bkap}(\Gt)$ such that 
\begin{subequations}\label{eqn:weak}
\begin{align}\label{eq:weak1}
&\bigl(\mathscr{V},~\varphi\bigr)_{\Gt}-\bigl(\nabs\varkappa,~\nabs\varphi\bigr)_{\Gt} + \bigl([\varkappa-\bkap]|\nabs\vec\nu|^2,~\varphi\bigr)_{\Gt}\nn\\
&\qquad\;-\frac{1}{2}\bigl([\varkappa-\bkap]^2\varkappa,~\varphi\bigr)_{\Gt}=0\qquad\forall\varphi\in  \bV_{0}(\Gt),\\[0.5em]
&\bigl(\partial_t^\circ\varkappa,~\chi\bigr)_{\Gt}+\frac{1}{2}\bigl(\nabs\cdot\mathscr{\vv V},~[\varkappa-\bkap]\,\chi\bigr)_{\Gt} -\frac{1}{2}\mathscr{A}_{\Gt}(\mathscr{\vv V},~\varkappa-\bkap,~\chi)\nn\\
&\qquad +\bigl(\nabs\mathscr{V},~\nabs\chi\bigr)_{\Gt} -\bigl(\mathscr{V}|\nabs\vec\nu|^2,~\chi\bigr)_{\Gamma(t)}\nn\\
&\qquad +\frac{1}{2}\bigl(\mathscr{V},~[\varkappa-\bkap]\,\varkappa\,\chi\bigr)_{\Gt}=0\qquad\forall\chi\in \bV_{\partial_1}(\Gt),\label{eq:weak2}\\[0.5em]
\label{eq:weak3}
&\bigl(\mathscr{\vv V}\cdot\vec \nu,~\xi\bigr)_{\Gt} - \bigl(\mathscr{V},~\xi\bigr)_{\Gt} = 0\qquad\forall\xi\in \bV_{0}(\Gt),\\[0.5em]
&\bigl(\kappa\,\vec\nu,~\vec\eta\bigr)_{\Gt}+\bigl(\nabs{\vec\id },~\nabs\vec\eta\bigr)_{\Gt} =0\qquad\forall\vec\eta\in [\bV_{0}(\Gt)]^d.\label{eq:weak4}
\end{align} 
\end{subequations}
Here \eqref{eq:weak1} is obtained by taking the inner product of \eqref{eq:rw1} with $\varphi\in \bV_0(\Gt)$ and then applying integration by parts.  Taking the inner product of \eqref{eq:rw2} with $\chi\in \bV_{\partial_1}(\Gt)$,  integrating by parts and recalling \eqref{eq:anti} as well as  $\vec\mu\cdot\nabs\mathscr{V}=0$ on $\partial_2\Gt$, we get \eqref{eq:weak2}. Equations \eqref{eq:weak3} and \eqref{eq:weak4} are due to \eqref{eq:rw3} and \eqref{eq:rw4}, respectively.  

\begin{thm} \label{thm:weakstab} A weak solution of \eqref{eqn:weak} 
satisfies the energy law 
\begin{equation}\label{eq:energylawweak}
\frac{\rm d}{\rm d t}E_{\bkap}(\varkappa(t), \vec x(t))+ \bigl(\mathscr{V},~\mathscr{V}\bigr)_{\Gt}=0.
\end{equation}
\end{thm}
\begin{proof}
It follows from \eqref{eq:transport} that
\begin{equation}
\frac{1}{2}\ddt\Bigl([\varkappa-\bkap]^2,~1\Bigr)_{\Gt} = \bigl(\partial_t^\circ\varkappa,~\varkappa-\bkap\bigr)_{\Gamma(t)} +\frac{1}{2}\Bigl(\nabs\cdot\mathscr{\vv V},~(\varkappa-\bkap)^2\Bigr)_{\Gt}.
\end{equation}
We next choose $\varphi = \mathscr{V}$ in \eqref{eq:weak1}, $\chi=\varkappa-\bkap$ in \eqref{eq:weak2}, and then combine the two equations to obtain that
\begin{equation}
\bigl(\mathscr{V},~\mathscr{V}\bigr)_{\Gt} + \bigl(\partial_t^\circ\varkappa,~\varkappa-\bkap\bigr)_{\Gt} +\frac{1}{2}\Bigl(\nabs\cdot\mathscr{\vv V},~(\varkappa-\bkap)^2\Bigr)_{\Gt}=0,
\end{equation}
which leads to \eqref{eq:energylawweak} straightforwardly. 
\end{proof}

We observe that Theorem \ref{thm:weakstab} holds for an arbitrary choice of $\varkappa(\cdot, 0)$. When it is chosen in a consistent manner, i.e., as the mean curvature of $\Gamma(0)$, we obtain that  \eqref{eqn:weak} describes the $L^2$-gradient flow of the energy \eqref{eq:Tenergy}. 

\section{Parametric finite element approximations}\label{sec:pfem}

 We employ a uniform partition of the time interval  as $[0,T]=\cup_{m=1}^M[t_{m-1}, t_m]$ with $t_m = m\ttau$ and $\ttau = \frac{T}{M}$.  Let $\Gm$ be a $(d-1)$-dimensional polyhedral surface to approximate the hypersurface $\Gamma(t_m)$ such that
\begin{equation}
 	\Gamma^m: = \cup_{j=1}^{J}\overline{\sigma_j^m}\quad\mbox{with}\quad \mT^m=\{\sigma_j^m\}_{j=1}^J,\quad \mQ^m=\{\vec q_k^m\}_{k=1}^{K},\label{eq:GammaD}
 \end{equation}
 where $\mT^m$ is a collection of mutually disjoint $(d-1)$-simplices in $\bR^d$,  and $\mQ^m$ is a collection of vertices of $\Gm$. Denote by $\partial_\alpha\Gm=\partial_\alpha\Gamma^0$ the approximation of the fixed boundary $\partial_\alpha\Gamma(0)$, $\alpha =1,2$.  In addition, let
$\left\{\vec q_{j_k}^{m}\right\}_{k=0}^{d-1}$ be the vertices of $\sigma_j^{m}$, ordered with the same orientation for all $\sigma\in \mT^m$.  
For simplicity, we denote $\sigma_j^{m}=\Delta\left\{\vec q_{j_k}^{m}\right\}_{k=0}^{d-1}$. Then we introduce the unit normal $\vec{\nu}^m$ to $\Gamma^m$; that is,
\begin{equation}\label{eq:vG}
\vec{\nu}^m_{j} := \vec{\nu}^m \mid_{\sigma^{m}_j} :=
\frac{\vec N\{\sigma_j^{m}\}}{
|\vec N\{\sigma_j^{m}\}|}\quad\mbox{ with}\quad \vec N\{\sigma_j^{m}\}=( \vec{q}^{m}_{j_1} - \vec{q}^{m}_{j_0} ) \wedge \ldots \wedge
( \vec{q}^{m}_{j_{d-1}} - \vec{q}^{m}_{j_0}),
\end{equation}
where $\wedge$ is the wedge product and $\vec N\{\sigma_j^{m}\}$ is the orientation vector of $\sigma_j^{m}$.  

On the polyhedral surface $\Gamma^m$, we define the finite element spaces
\begin{align}
	\bV^h(\Gm) &:=\bigl\{\chi\in C(\Gamma^m): \chi|_{\sigma}\;\;\mbox{is affine}\quad\forall  \sigma\in\mT^m \bigr\},\nn\\
	\bV_{\partial_\alpha}^h(\Gm)&:=\bV_{\partial_\alpha}(\Gm)\cap\bV^h(\Gm),\quad \alpha=1,2;\nn\\
	\bV^h_0(\Gm) &:=\bV_0(\Gm)\cap\bV^h(\Gm); \nn
\end{align}
and the sets
\begin{align}
	\bK_{\partial_\alpha}^h(\Gm)&=\bK_{\partial_\alpha}(\Gm)\cap\bV^h(\Gm), \quad \alpha=1,2;\nn\\
		\bK_{\bkap }^h(\Gm)&:=\bK_{\bkap}(\Gm)\cap\bV^h(\Gm). \nn 
\end{align}
To approximate the inner product $\ipd{\cdot,\cdot}_{\Gm}$, we introduce the mass-lumped approximation over
the current polyhedral surface $\Gamma^m$ via 
\begin{equation}
\ipd{u, v}_{\Gm}^h:=  
\frac{1}{d}\sum_{j=1}^{J} \mH^{d-1}(\sigma^{m}_j)
\sum_{k=0}^{d-1} 
\underset{\sigma^{m}_j\ni \vec{p}\to \vec{q}^{m}_{j_k}}{\lim}\, 
(u\cdot v)(
\vec{p}),\label{eq:tprule}
\end{equation}
where $u,v$ are piecewise continuous, with possible jumps
across the edges of $\sigma\in\mT^m$, and 
$\mH^{d-1}(\sigma^{m}_j)= \frac{1}{(d-1)!}\,|\vec N\{\sigma_j^{m}\}|$ 
is the measure of $\sigma^{m}_j$.

The discrete normals $\vec\nu^m$ in \eqref{eq:vG} are piecewise constant. For later use, we follow the work in \cite{BGN08parametric, BGN08willmore} and introduce the vertex normal vector $\vec\omega^{m}\in [\bV^h(\Gm)]^d$.  
That is, we let $\vec\omega^m\in  [\bV^h(\Gm)]^d$ be the the mass-lumped $L^2$--projection of $\vec\nu^m$ onto $[\bV^h(\Gm)]^d$, i.e., 
\begin{equation} \label{eq:nuhomegah}
\ipd{\vec\omega^m, \vec\eta^h}_{\Gm}^h = \ipd{\vec\nu^m,~\vec\eta^h}_{\Gm}^h = \ipd{\vec\nu^m,~\vec\eta^h}_{\Gm}\quad\forall\vec\eta^h\in [\bV^h(\Gm)]^d.
\end{equation}
It is not difficult to show that the following identity holds
\begin{equation*}
\ipd{\chi\,\vec\omega^m,~\vec\eta^h}_{\Gm}^h = \ipd{\chi\,\vec\nu^m,~\vec\eta^h}_{\Gm}^h\qquad\forall\chi\in \bV^h(\Gm),\quad\vec\eta^h\in [\bV^h(\Gm)]^d.
\end{equation*}

\subsection{The linear  method} 
Inspired by the ALE framework \cite{GNZ24ALE},  we first define $\mathscr{\vv V}^m\in[\bV^h(\Gm)]^d$ as the discrete vertex velocity on $\Gamma^m$
\begin{equation} \label{eq:Vm}
\mathscr{\vv V}^m(\vec q_k^m) = \frac{\vec q_k^{m} - \vec q_k^{m-1}}{\ttau},\qquad 1\leq k\leq K,
\end{equation}
and then introduce the discrete mapping $\vec\Phi^m\in[\bV^h(\Gm)]^d$ as
\begin{equation}
\vec\Phi^m = \vec\id - \ttau\,\mathscr{\vv V}^m,\quad\mbox{so that}\quad \vec\Phi^m(\vec q_k^m) = \vec q_k^{m-1},\quad \vec\Phi^m(\Gamma^m) = \Gamma^{m-1}.\label{eq:mapping}
\end{equation}

Let $\widehat{\mathcal{L}}^m(\vec z)=
[\nabs\vec\Phi^m(\vec z)]^T\nabs\vec\Phi^m(\vec z) \in \bR^{d\times d}$ 
for $\vec z\in\Gm$. For a fixed $\vec z\in\Gm$,
this induces a linear operator 
$\mathcal{L}^m(\vec z): T_{\vec z}\Gm\mapsto T_{\vec z}\Gm$ on the tangent 
space $T_{\vec z}\Gm$, such that 
$[\mathcal{L}^m(\vec z)]\,\vec t = [\widehat{\mathcal{L}}^m(\vec z)]\,\vec t$
for all $\vec t \in T_{\vec z}\Gm$. 
The surface Jacobian determinant of the mapping $\vec\Phi^m$ is then given by
\begin{equation}
	\mathcal{J}^m=\sqrt{{\rm det}(\mathcal{L}^m)} ,
\label{eq:jmdef}
\end{equation}
which describes local area change and is a constant on each $\sigma\in\mT^m$. 
In particular,
\begin{equation} \label{eq:jmint}
\int_{\Gamma^{m-1}} f\, \dH^{d-1} = \int_{\Gamma^m} f \circ \vec\Phi^m\,
\mathcal{J}^m \, \dH^{d-1} \qquad\forall f\in L^1(\Gamma^{m-1}).
\end{equation}
We observe that \eqref{eq:Vm} is a consistent temporal approximation
of $\mathscr{\vv V}$. Hence it is natural to assume that
$\nabs\cdot \mathscr{\vv V}^m$ remains bounded throughout the evolution. 
This yields the following lemma
for $\mathcal{J}^m$. 
\begin{lemma} \label{lem:Jm}
	If $\max_{1\leq m\leq M}|\nabs\mathscr{\vv V}^m|\leq C$ for some constant $C$ independent of the time step $\ttau$, then it holds that
\begin{equation}\label{eq:Jmexp}
	\sqrt{\mathcal{J}^m} = 1 - \frac{1}{2}\ttau\,\nabs\cdot \mathscr{\vv V}^m + O(\ttau^2), \quad m = 1,\cdots, M.
\end{equation}
\end{lemma}	
\begin{proof}
We restrict our attention to a simplex $\sigma\in\mathscr{T}^m$, and let $\{\vec t_1,\cdots, \vec t_{d-1}\}$ be an orthonormal basis of the tangential space of  $\sigma$. In order to compute ${\rm det}(\mathcal{L}^m)$, we consider the matrix of $\mathcal{L}^m$ with respect to the  basis $\{\vec t_1,  \cdots, \vec t_{d-1}\}$ and denote it by $\mat{M}^m = (M_{ij})_{i,j=1,\cdots,d-1}$. Then it holds that
\begin{equation}
\mathcal{L}^m\,\vec t_j = \sum_{i=1}^{d-1} M_{ij}\vec t_i\qquad j = 1,\cdots, d-1,
\end{equation}
which by the orthonormality  implies that
\[M_{ij} = \vec t_i^T(\mathcal{L}^m\,\vec t_j) =[\nabs\vec\Phi^m\,\vec t_i]^T\nabs\vec\Phi^m\,\vec t_j ,\quad i,j = 1,\cdots, d-1.\] 

Recalling \eqref{eq:mapping} we can recast the matrix $\mat{M}^m$ as
\begin{align}
	\mat{M}^m &= \left(\{\vec t_i -\ttau[\nabs\mathscr{\vv V}^m]\vec t_i\}^T\{\vec t_j -\ttau[\nabs\mathscr{\vv V}^m]\vec t_j\}\right)_{i,j=1,\cdots, d-1}\nn\\
	& = \left(\delta_{ij} - \ttau\,\left\{\vec t_j^T[\nabs\mathscr{\vv V}^m]\vec t_i +\vec t_i^T[\nabs\mathscr{\vv V}^m]\vec t_j\right\} + O(\ttau^2)\right)_{i,j=1,\cdots, d-1}\nn\\
	& =\mat{I}_{d-1}- \ttau\,\mat{A}^m + O(\ttau^2),
	\label{eq:jm3}
\end{align}
where $\mat{I}_{d-1}\in\bR^{(d-1)\times(d-1)}$ is the identity  matrix, and  
\[\mat{A}^m = \left(\vec t_j^T[\nabs\mathscr{\vv V}^m]\vec t_i +\vec t_i^T[\nabs\mathscr{\vv V}^m]\vec t_j\right)_{i,j=1,\cdots, d-1}.\]
Now using \eqref{eq:jm3} yields,  on recalling
$\nabs\cdot\mathscr{\vv V}^m = \sum_{i=1}^{d-1} \vec t_i \cdot \partial_{\vec t_i} \mathscr{\vv V}^m$, that
\begin{align*}
	\left(\mathcal{J}^m\right)^2 &={\rm det}(\mathcal{L}^m)={\rm det}(\mat{M}^m)={\rm det}(\mat{I}_{d-1} - \ttau\,\mat{A}^m + O(\ttau^2))\nn\\
	& = 1 -\ttau\,{\rm tr}(\mat{A}^m) + O(\ttau^2)\nn\\
	& = 1- 2\ttau\sum_{i=1}^{d-1}\vec t_i^T[\nabs\mathscr{\vv V}^m]\vec t_i + O(\ttau^2) = 1 - 2\ttau\,\nabs\cdot \mathscr{\vv V}^m + O(\ttau^2),
\end{align*}
where for the equality at the second line we have observed that
${\rm det}(\mat{I}_{2}+\delta\mat{B}) = 1 + \delta\,{\rm tr}(\mat{B}) 
+ \delta^2 {\rm det}(\mat{B})$ for any $\mat{B}\in\bR^{2\times2}$. Using the expansion of $(1-\delta x)^{\tfrac{1}{4}}= 1 - \tfrac{1}{4} \delta x + O((\delta x)^2)$,  we get \eqref{eq:Jmexp} as claimed.
\end{proof}

\vspace{0.5cm}

Denote by $\varkappa^m, \mathscr{V}^m$ and $\kappa^m$ the numerical approximations of $\varkappa(\cdot, t_m), \mathscr{V}(\cdot, t_m)$ and $\kappa(\cdot, t_m)$, respectively.  
We are now ready to propose our novel finite element scheme for the weak formulation \eqref{eqn:weak} as follows.   Given the initial non-degenerate polyhedral surface $\Gamma^0$ and $\varkappa^0_{\Gamma^0}:=\varkappa^0\in \bK_{\bkap}^h(\Gamma^0)$, we set $\Gamma^{-1}=\Gamma^0$ with $\mathcal{J}^0=1$. Then for each $m\geq 0$, we first obtain an appropriate approximation $\mathcal{W}^m \in L^\infty(\Gamma^m)$ of $|\nabs\vec\nu|^2$ at $t_m$. Precisely, we let
\begin{equation} \label{eq:Wm}
\mathcal{W}^m =| \nabs\vec v^m|^2\quad\mbox{with}\quad  \vec v^m(\vec q) =\frac{\vec\omega^m(\vec q)}{|\vec\omega^m(\vec q)|}\qquad\forall\vec q\in \mQ^m, 
\end{equation}
where $\vec\omega^m$ is the vertex normal defined in \eqref{eq:nuhomegah}, which from now on we assume to be nondegenerate, and $\vec v^m\in[\bV^h(\Gm)]^d$ is the normalized vertex normal. We recall from \cite[Remark 65]{Barrett20} that $\vec\omega^m$ can only be degenerate in some very pathological
cases. For example, for surfaces without
self-intersections it will always be
nondegenerate. The next step of the scheme consists of two stages: computing the normal velocity and applying the normal velocity together with a tangential velocity. 

  \vspace{0.6em}
\noindent {\bf  Stage 1}:  Find $\left(\mathscr{V}^{m+1}, \varkappa^{m+1}\right) \in \bV_{0}^h(\Gamma^m)\times \bK^h_{\partial_{1}}(\Gamma^m)$ such that
\begin{subequations}\label{eqn:fdVn}
\begin{align}\label{eq:fd1}
&\ipd{\mathscr{V}^{m+1},~\varphi^h}_{{\Gm}}-\bigl(\nabs\varkappa^{m+1},~\nabs\varphi^h\bigr)_{\Gm} + \bigl(\mathcal{W}^m\,[\varkappa^{m+1}-\bkap],~\varphi^h\bigr)_{\Gm}\nn\\
&\qquad-\frac{1}{2}\bigl([\varkappa^m_{{\Gamma^m}}-\bkap]\,\varkappa_{\Gamma^m}^m\,[\varkappa^{m+1}-\bkap],~\varphi^h\bigr)_{\Gm}=0,\\[0.5em]
&\Bigl(\frac{\varkappa^{m+1}-\bkap - (\varkappa_{\Gamma^m}^m-\bkap)\,\sqrt{\mathcal{J}^m}}{\ttau},~\chi^h\Bigr)_{\Gm}-\frac{1}{2}\mathscr{A}_{\Gm}(\mathscr{\vv V}^m,~\varkappa^{m+1}-\bkap,~\chi^h)\nn\\
&\qquad +\bigl(\nabs\mathscr{V}^{m+1},~\nabs\chi^h\bigr)_{\Gm} -\bigl(\mathcal{W}^m\,\mathscr{V}^{m+1},~\chi^h\bigr)_{\Gm}\nn\\
&\qquad+\frac{1}{2}\bigl(\mathscr{V}^{m+1},~[\varkappa^m_{\Gamma^m}-\bkap]\,\varkappa_{\Gamma^m}^m\,\chi^h\bigr)_{\Gm}=0,\,\label{eq:fd2}
\end{align} 
\end{subequations}
for $(\varphi^h, \chi^h)\in \bV_{0}^h(\Gm)\times\bV_{\partial_{1}}^h(\Gm)$.  

  \vspace{0.6em}
\noindent{\bf Stage 2}: Find $\kappa^{m+1}\in \bK^h_{\bkap}(\Gm)$ and $\vec X^{m+1}\in [\bV^h(\Gm)]^d$  with $(\vec X^{m+1}-\vec\id)\in [\bV^h_0(\Gm)]^d$ such that
\begin{subequations}\label{eqn:fdVt}
\begin{align}
	\label{eq:fd3}
&\Bigl(\frac{\vec X^{m+1}-\vec\id }{\ttau}\cdot\vec \nu^m,~\xi^h\Bigr)_{\Gm}^h - \bigl(\mathscr{V}^{m+1},~\xi^h\bigr)_{\Gm} = 0,\\[0.5em]
&\Bigl(\kappa^{m+1}\,\vec\nu^m,~\vec\eta^h\Bigr)_{\Gm}^h+\Bigl(\nabs{\vec X^{m+1}},~\nabs\vec\eta^h\Bigr)_{\Gm}= 0,\label{eq:fd4}
\end{align}
\end{subequations}
$(\xi^h,~\vec\eta^h)\in\bV_0^h(\Gm)\times [\bV^h_0(\Gm)]^d$, and then set $\Gamma^{m+1}=\vec X^{m+1}(\Gm)$, as well as $\varkappa_{\Gamma^{m+1}}^{m+1} = \varkappa^{m+1}\circ\vec\Phi^{m+1}\in \bK_{\partial_1}^h(\Gamma^{m+1})$. Observe that we employ mass-lumping in 
\eqref{eqn:fdVt} in order to obtain the desired BGN tangential motion
\cite{Barrett20}. Other choices for the tangential motion can easily be
incorporated, see \S\ref{sec:fi} below. 

 \vspace{0.6em}
 \begin{rem}\label{rem:tp}
By \eqref{eq:Jmexp}, it is expected that
\begin{align}
&\Bigl(\frac{\varkappa^{m+1}-\bkap-(\varkappa_{\Gamma^m}^m-\bkap)\,\sqrt{\mathcal{J}^m}}{\ttau},~\chi^h\Bigr)_{\Gamma^m} \nn\\
&=\ipd{\frac{\varkappa^{m+1}-\varkappa_{\Gamma^m}^m}{\ttau},~\chi^h}_{\Gamma^m} +\frac{1}{2}\ipd{\nabla_s\cdot\mathscr{\vv V}^m,~[\varkappa^m_{\Gamma^m}-\bkap],\chi^h}_{\Gamma^m} +O(\ttau), \label{eq:ALEapp}
\end{align}
which shows that the first term of \eqref{eq:fd2} is a consistent temporal discretization of the first two terms in \eqref{eq:weak2}. 
We remark that it does not seem possible to rigorously prove the
assumption on the uniform boundedness of $\nabla_s \cdot \mathscr{\vv V}^m$
that is needed in Lemma~\ref{lem:Jm}. However, in practice the scheme
\eqref{eqn:fdVn}, \eqref{eqn:fdVt} works very well. Moreover, in 
\S\ref{sec:fi} we briefly discuss a nonlinear scheme that is also 
unconditionally energy stable, and for which Lemma~\ref{lem:Jm}
is not needed.
\end{rem}

We note that the introduced scheme \eqref{eqn:fdVn}, \eqref{eqn:fdVt} leads to a decoupled system of two linear systems of equations. In the case of closed planar curves with $\bkap=0$, the scheme collapses to the linear scheme that was recently introduced by the present authors in \cite{GNZ25willmore}, if we replace \eqref{eq:Wm} with $\mathcal{W}^m = (\varkappa^m)^2$ .

\subsection{The properties of the method}

We next aim to show that the two linear systems from \eqref{eqn:fdVn}, \eqref{eqn:fdVt} admit a unique solution. We denote by
\begin{align}
	\mQ_{\circ}^m:=\left\{\vec q\in \mQ^m:\; \vec q\notin\partial\Gm \right\}, 
\end{align}
  the set of the interior vertices of $\Gamma^m$. Then we have the following theorem. 
\begin{thm} [existence and uniqueness]  
Assume that $\mH^{d-1}(\sigma) > 0$ for all $\sigma \in \mT^m$, and that
$|\vec\omega^m(\vec q)| > 0$ for all $\vec q \in \mQ^m$.
Then the linear system \eqref{eqn:fdVn} admits a unique solution $(\mathscr{V}^{m+1}, ~\varkappa^{m+1})$.  Moreover, on assuming that
\begin{enumerate}[label=$(\mathbf{A \arabic*})$, ref=$\mathbf{A \arabic*}$] 
\item \label{assupI} $\partial\Gm \not=\emptyset$ or ${\rm dim\; span}\left(
\bigl\{\vec\omega^{m}(\vec q)\bigr\}_{\vec q\in \mQ_{\circ}^m}\right)=d$,
\end{enumerate}
then the linear system \eqref{eqn:fdVt} admits a unique solution for any given $\mathscr{V}^{m+1}$. 
\end{thm}

\begin{proof}
We first consider the linear system \eqref{eqn:fdVn}. Since  the number of unknowns equals the number
of equations, it suffices to show that the corresponding homogeneous system has only the zero solution. Therefore, we consider the following homogeneous system for $(\mathscr{V}, \varkappa)\in \bV^h_{0}(\Gm)\times \bV^h_{\partial_1}(\Gm)$ such that
\begin{subequations}\label{eqn:homo}
\begin{align}\label{eq:homo1}
&\ipd{\mathscr{V},~\varphi^h}_{{\Gm}}-\bigl(\nabs\varkappa,~\nabs\varphi^h\bigr)_{\Gm} + \bigl(\mathcal{W}^m\,\varkappa,~\varphi^h\bigr)_{\Gm}\nn \\ &\qquad\qquad\;
-\frac{1}{2}\bigl([\varkappa_{\Gamma^m}^m-\bkap]\,\varkappa_{\Gamma^m}^m\,\varkappa,~\varphi^h\bigr)_{\Gm}=0,\\[0.5em]
&\Bigl(\frac{\varkappa}{\ttau},~\chi^h\Bigr)_{\Gm} +\bigl(\nabs\mathscr{V},~\nabs\chi^h\bigr)_{\Gm} -\bigl(\mathcal{W}^m\,\mathscr{V},~\chi^h\bigr)_{\Gm}\nn\\
&\qquad\quad  -\frac{1}{2}\mathscr{A}_{\Gm}(\mathscr{\vv V}^m,~\varkappa,~\chi^h) +\frac{1}{2}\bigl(\mathscr{V},~[\varkappa_{\Gamma^m}^m-\bkap]\,\varkappa_{\Gamma^m}^m\,\chi^h\bigr)_{\Gm}=0,\,\label{eq:homo2}
\end{align}
\end{subequations}
for $(\varphi^h, \chi^h)\in\bV^h_{0}(\Gm)\times\bV^h_{\partial_{1}}(\Gm)$. We then set $\varphi^h =\mathscr{V}$ in \eqref{eq:homo1} and $\chi^h =\varkappa$ in \eqref{eq:homo2} to obtain, on recalling \eqref{eq:antisym}, that
\begin{equation}
	\ttau \ipd{\mathscr{V}, ~\mathscr{V}}_{\Gm} + \ipd{\varkappa,~\varkappa}_{\Gm}=0,
\end{equation}
which immediately implies that $(\mathscr{V}, ~\varkappa)=(0, 0)$. Thus \eqref{eqn:fdVn} has a unique solution. 

The proof of the unique solvability of the linear system \eqref{eqn:fdVt} in
the case $\partial\Gamma^m=\emptyset$ is given in \cite[Lemma~66]{Barrett20}.
In the case $\partial\Gamma^m\not=\emptyset$ the proof is even simpler, and so
we omit it here.
\end{proof}

We next show that the scheme \eqref{eqn:fdVn}, \eqref{eqn:fdVt} satisfies a stability bound that mimics the energy dissipation law in \eqref{eq:energylaw} on the fully discrete level. 

\begin{thm}[unconditional stability] \label{thm:stability} Let $(\mathscr{V}^{m+1}, \varkappa^{m+1}, \mathscr{\vv V}^{m+1}, \kappa^{m+1})$ be a solution to \eqref{eqn:fdVn}, \eqref{eqn:fdVt}. Then for $m=0,1,\cdots, M-1$, it holds that
\begin{equation}\label{eq:dES}
\frac{1}{2}\norm{\varkappa^{m+1}-\bkap}_{m}^2 + \ttau\norm{\mathscr{V}^{m+1}}_m^2 \leq \frac{1}{2}\norm{\varkappa^{m}-\bkap}_{m-1}^2,
\end{equation}
where  $\norm{\cdot}_m$ is the norm induced by the inner product $(\cdot,\cdot)_{_\Gm}$.
\end{thm}

\begin{proof}
We set $\varphi^h = \ttau\,\mathscr{V}^{m+1}$ in \eqref{eq:fd1}, $\chi^h = \ttau\,(\varkappa^{m+1}-\bkap)$ in \eqref{eq:fd2} and combine these equations to obtain that
\begin{align}
&\ttau\bigl(\mathscr{V}^{m+1},~\mathscr{V}^{m+1}\bigr)_{\Gm} + \bigl(\varkappa^{m+1}-\bkap-(\varkappa_{\Gamma^m}^m - \bkap)\,\sqrt{\mathcal{J}^m},~\varkappa^{m+1}-\bkap\bigr)_{\Gm}=0.\label{eq:es1}
\end{align}

Using the inequality $a(a-b)\geq \frac{1}{2}(a^2-b^2)$ with $a=\varkappa^{m+1}-\bkap$ and $b = (\varkappa^m_{\Gamma^m}-\bkap)\sqrt{\mathcal{J}^m}$, we have,
on recalling \eqref{eq:jmint}, that 
\begin{align}
 &\bigl(\varkappa^{m+1}-\bkap-(\varkappa_{\Gamma^m}^m - \bkap)\,\sqrt{\mathcal{J}^m},~\varkappa^{m+1}-\bkap\bigr)_{\Gm}\nn\\ & \qquad
\geq \frac{1}{2}\ipd{[\varkappa^{m+1}-\bkap]^2, ~1}_{\Gamma^m} - \ipd{[\varkappa^{m}\circ\vec\Phi^m-\bkap]^2,~\mathcal{J}^m}_{\Gm}\nn\\ & \qquad
 = \frac{1}{2}\norm{\varkappa^{m+1}-\bkap}_m^2 - \frac{1}{2}\norm{\varkappa^{m}-\bkap}_{m-1}^2.\label{eq:es2}
 \end{align}
Inserting \eqref{eq:es2} into \eqref{eq:es1} yields \eqref{eq:dES} as claimed. 
\end{proof}

Note that the stability results in both \eqref{eq:energylawweak} and \eqref{eq:dES} only control the normal component of the velocity rather than the full velocity. This is consistent with the fact that the gradient flow structure in (2.6) is entirely determined by the normal velocity $\mathscr{V}$ on $\Gt$.

\subsection{Further insights} \label{sec:fi}
Building upon the framework for the introduced scheme \eqref{eqn:fdVn}, \eqref{eqn:fdVt}, several additional observations and considerations arise that illuminate the strengths and potential extensions of the proposed method.

First of all, alternative choices for $\mathcal{W}^m$ to \eqref{eq:Wm} can
easily be used, and do not affect any of the theoretical properties shown for
the scheme. For example, for $d=2$ one can use 
$\mathcal{W}^m = (\varkappa^m)^2$.

Secondly, in the stability estimate \eqref{eq:dES}, we notice that the discrete energy is defined by integrating $\varkappa^{m+1}$ over the surface $\Gm$. 
Upon modifying the introduced scheme to a nonlinear method, it is possible to
derive a stability estimate where the discrete energy is instead given by 
$\frac{1}{2}\norm{\varkappa^{m+1}_{\Gamma^{m+1}}-\bkap}_{m+1}^2$, with $\varkappa^{m+1}_{\Gamma^{m+1}} = \varkappa^{m+1}\circ\vec\Phi^{m+1}\in \bK_{\partial_1}^h(\Gamma^{m+1})$. 
On recalling $\nabs\cdot\mathscr{\vv V} = \nabs\vec\id: \nabs\mathscr{\vv V}$,
the new scheme would approximate \eqref{eq:weak2} directly via
\begin{align}
&\quad\ipd{\frac{\varkappa^{m+1} - \varkappa_{\Gamma^m}^m}{\ttau}, \chi^h}_{\Gm} + \frac{1}{2}\ipd{\nabs\vec X^{m+1}:\nabs\left[\frac{\vec X^{m+1}-\vec\id}{\ttau}\right],[\varkappa^{m+1}-\bkap]\chi^h}_{\Gamma^m}\nn\\
&\qquad\quad  -\frac{1}{2}\mathscr{A}_{\Gm}(\mathscr{\vv V}^m, \varkappa^{m+1}-\bkap,~\chi^h) +\bigl(\nabs\mathscr{V}^{m+1},~\nabs\chi^h\bigr)_{\Gm}
 \nn\\
  &\qquad\qquad -\bigl(\mathcal{W}^m\,\mathscr{V}^{m+1},~\chi^h\bigr)_{\Gm} +\frac{1}{2}\bigl(\mathscr{V}^{m+1},~[\varkappa_{\Gamma^m}^m-\bkap]\,\varkappa_{\Gamma^m}^m\,\chi^h\bigr)_{\Gm}=0\,,\label{eq:nnfd2}	
\end{align}
rather than using \eqref{eq:fd2} in combination with
\eqref{eq:ALEapp} and Lemma~\ref{lem:Jm}. 
Replacing \eqref{eq:fd2} with \eqref{eq:nnfd2} leads to a nonlinear scheme, which couples {\bf Stage 1} and {\bf Stage 2} in view of the presence of $\vec X^{m+1}$ in \eqref{eq:nnfd2}.  Nevertheless, this new method satisfies the stability estimate
\begin{equation}\label{eq:dESn}
\frac{1}{2}\norm{\varkappa_{\Gamma^{m+1}}^{m+1}-\bkap}_{m+1}^2 + \ttau\norm{\mathscr{V}^{m+1}}_m^2 \leq \frac{1}{2}\norm{\varkappa_{\Gamma^m}^{m}-\bkap}_{m}^2,
\end{equation}
which is slightly different from that in \eqref{eq:dES}. The proof of \eqref{eq:dESn} is similar to that of \eqref{eq:dES}, except that instead of \eqref{eq:es1}, we get  
\begin{align}
\ttau\ipd{\mathscr{V}^{m+1}, \mathscr{V}^{m+1}}_{\Gm} + \ipd{\varkappa^{m+1}-\varkappa^m_{\Gamma^m},~\varkappa^{m+1}-\bkap}_{\Gm} \nn\\
+ \frac{1}{2}\ipd{\nabs\vec X^{m+1}: \nabs[\vec X^{m+1}-\vec\id], [\varkappa^{m+1}-\bkap]^2}_{\Gm}
 = 0.\label{eq:nstab}
\end{align}
Using the inequality $a(a-b)\geq \frac{1}{2}(a^2-b^2)$ with $a=\varkappa^{m+1}-\bkap$ and $b=\varkappa^m_{\Gamma^m}-\bkap$ we obtain
\begin{align}
&\ipd{\varkappa^{m+1}-\varkappa^m_{\Gamma^m}, \varkappa^{m+1}-\bkap}_{\Gm}\geq \frac{1}{2}\norm{\varkappa^{m+1}-\bkap}_{m}^2 - \frac{1}{2}\norm{\varkappa^m_{\Gamma^m}-\bkap}_m^2.\label{eq:nstab1}
\end{align}
In addition, it follows from Appendix~\ref{sec:nstab2} that
\begin{equation}
\ipd{\nabs\vec X^{m+1}: \nabs[\vec X^{m+1}-\vec\id], [\varkappa^{m+1}-\bkap]^2}_{\Gm}
\geq \norm{\varkappa^{m+1}_{\Gamma^{m+1}}-\bkap}_{m+1}^2 - \norm{\varkappa^{m+1} - \bkap}_{m}^2.\label{eq:nstab2}
\end{equation}
Combining \eqref{eq:nstab}, \eqref{eq:nstab1} and \eqref{eq:nstab2} yields the desired result \eqref{eq:dESn}.

Finally, the proof of Theorem~\ref{thm:stability} shows the flexibility of our approach in choosing the tangential motion in {\bf Stage 2}, while still preserving the energy decay property on the discrete level. Therefore, it is possible to replace \eqref{eqn:fdVt} by other approaches, that either directly or indirectly fix the tangential velocities. An example for the latter is the minimal deformation rate (MDR) method from \cite{Hu22evolving}. To use the MDR method, we replace \eqref{eqn:fdVt} with
	\begin{subequations}\label{eqn:fdVtmdr}
\begin{align}
	\label{eq:fd3mdr}
&\Bigl(\frac{\vec X^{m+1}-\vec\id }{\ttau}\cdot\vec \nu^m,~\xi^h\Bigr)_{\Gm}^h - \bigl(\mathscr{V}^{m+1},~\xi^h\bigr)_{\Gm} = 0,\\[0.5em]
&\Bigl(\kappa^{m+1}\,\vec\nu^m,~\vec\eta^h\Bigr)_{\Gm}^h+\Bigl(\nabs\left[\frac{\vec X^{m+1}-\vec\id}{\ttau}\right],~\nabs\vec\eta^h\Bigr)_{\Gm}= 0,\label{eq:fd4mdr}
\end{align}
\end{subequations}
$(\xi^h,~\vec\eta^h)\in\bV_0^h(\Gm)\times [\bV^h_0(\Gm)]^d$.
This results in a different tangential velocity. Moreover, the variable 
$\kappa^{m+1}$ no longer approximates the mean curvature, see \cite{Hu22evolving} for further details.

\section{Numerical results}
\label{sec:num}

In this section, we present a variety of numerical examples for our introduced scheme \eqref{eqn:fdVn}, \eqref{eqn:fdVt}. 
We implemented the scheme within the
finite element toolbox ALBERTA, see \cite{Alberta}. The 
linear systems of equations arising at each time level
are solved with the help of the sparse factorization package UMFPACK, 
\cite{Davis04}.

Throughout the experiments,  we first construct an initial polyhedral surface $\Gamma_Y$ with $\vec Y=\vec\id|_{\Gamma_Y}\in [\bV^h(\Gamma_Y)]^d$ being the identity function.  To start our computations, we require the initial data $\varkappa^0$. In the case when the initial surface is part of a sphere of radius $r_0$, we define $\varkappa^0 = -\frac{d-1}{r_0}$.  Otherwise, we solve for the discrete mean curvature using the BGN method with a zero normal velocity and fixed boundary. That is, we find $(\delta\vec  Y^0, \kappa^0)\in [\bV_0^h(\Gamma_Y)]^d\times\bK_{\bkap}^h(\Gamma_Y)$ such that
\begin{subequations}\label{eqn:ic}
	\begin{align}
		&\ipd{\delta \vec Y^0\cdot\vec\nu_Y, \xi^h}_{\Gamma_Y}^h =0\qquad\forall\xi^h\in \bV_{0}^h(\Gamma_Y),\\
		&\ipd{\kappa^0\,\vec\nu_Y,~\vec\eta^h}_{\Gamma_Y}^h + \ipd{\nabs(\vec\id+\delta\vec  Y^0),~\nabs\vec\eta^h}_{\Gamma_Y}=0\qquad\forall\vec\eta^h\in [\bV_{0}^h(\Gamma_Y)]^d,
	\end{align}
	\end{subequations}
	where $\vec\nu_Y$ is the unit normal of $\Gamma_Y$ which follows \eqref{eq:vG} similarly. We next set $\vec X^0 = \vec\id|_{\Gamma_Y} + \delta\vec Y^0$ with $\Gamma^0 = \vec X^0(\Gamma_Y)$ and $\varkappa^0=\kappa^0$ as the required initial data.  For later use,  we introduce the time interpolation of the polyhedral surfaces via
	\begin{equation}
		\vec X_{h,\ttau} (\vec q, t) = \frac{t_{m+1}-t}{\ttau}\vec q + \frac{t-t_m}{\ttau}\vec X^{m+1}(\vec q),\quad\forall\vec q\in\mQ^m, \quad t\in[t_m, t_{m+1}].\label{eq:linertime}
	\end{equation}

\subsection{2d results} 
In this subsection, we will perform a convergence test for our scheme by considering the evolution of both closed and open planar curves in $\bR^2$.  The effects of different boundary conditions on the evolution of open curves will also be investigated. 

\begin{figure}[t]
	\centering
	\includegraphics[width=0.9\textwidth]{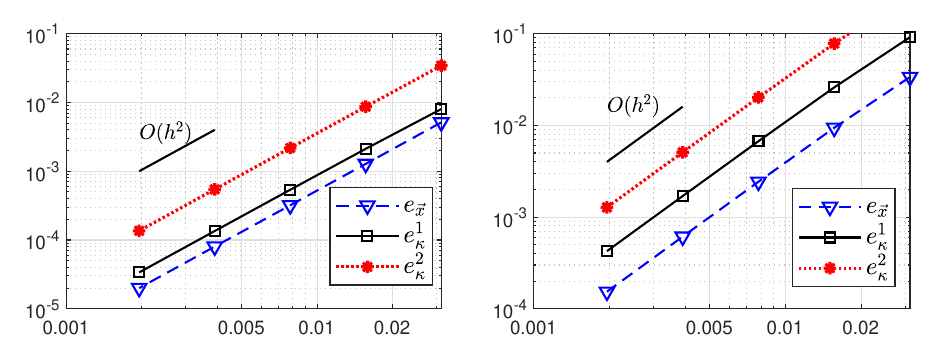}
	\caption{  Numerical errors for an expanding/shrinking circle with $\bkap=-0.5$ (left panel) and $\bkap=-2$ (right panel),  where $\ttau=(\frac{2^{5}\,h}{5})^2 $.}
	\label{fig:closederr}
\end{figure}

\vspace{0.3cm}
\noindent
{\bf   Example 1}: We begin with a convergence test for an initial circle of radius $r(0)=r_0$.  We note a closed circle with  radius $r(t)$ satisfying 
\begin{equation}
	r^\prime(t) = \frac{1}{2\,r(t)}\left(\frac{1}{r(t)}-\bkap\right)\left(\frac{1}{r(t)}+\bkap\right),\quad r(0) = r_0>0, \nonumber 
\end{equation}
is an exact solution of \eqref{eq:Willmore}. In the case of $\bkap\neq 0$ and $|\bkap\,r_0|\neq 1$,  the solution of the ordinary differential equation satisfies $\bkap^4 t+\bkap^2(r^2(t)-r_0^2) + \ln\left(\frac{1-\bkap^2\,r^2(t)}{1-\bkap^2\,r_0^2}\right)=0$. We choose a unit circle with $r_0= 1$ and consider  $\bkap = -0.5$ and $\bkap=-2$, meaning that the circle will expand/shrink towards a circle of radius 2 and 0.5, respectively. Initially  $\Gamma(0)$ is approximated by discretizing 
\begin{equation*}
\vec x(\rho,t)= \left(\begin{matrix}
\cos g(\rho)\\
\sin g(\rho)
\end{matrix}\right), \quad g(\rho) = -2\pi\rho + 0.1\sin(-2\pi\,\rho), \quad\rho\in\bR/\mathbb{Z},
\end{equation*}
uniformly according to $\rho$. Observe that $g(\rho)$ leads to a nonuniform initial distribution of nodes.  We introduce the errors over the time interval $[0,2]$ as
\begin{subequations}\label{eqn:dderror}
\begin{align}
\label{eq:Xerror}
e_{\vec x} &= \max_{1\leq m\leq M}\max_{\vec q\in\mQ^m}\left||\vec X^m(\vec q)| - r(t_m)\right|,\\
\label{eq:kaerror}
e_\kappa^1 &= \max_{1\leq m\leq M}\max_{\vec q\in\mQ^m}\left|\varkappa^m(\vec q) - \varkappa(t_m)\right|,\\
e_\kappa^2 &= \max_{1\leq m\leq M}\max_{\vec q\in\mQ^m}\left|\kappa^m(\vec q) - \varkappa(t_m)\right|.
\end{align}
\end{subequations}
The numerical results are reported in Fig. \ref{fig:closederr}, where we fix $\ttau=(\frac{2^{5}\,h}{5})^2 $ with $h=\frac{1}{J}$.  Here we observe a second-order convergence for our introduced scheme. 

\begin{figure}[!htp]
\vspace{-0.3cm}
	\centering
	\includegraphics[width=0.9\textwidth]{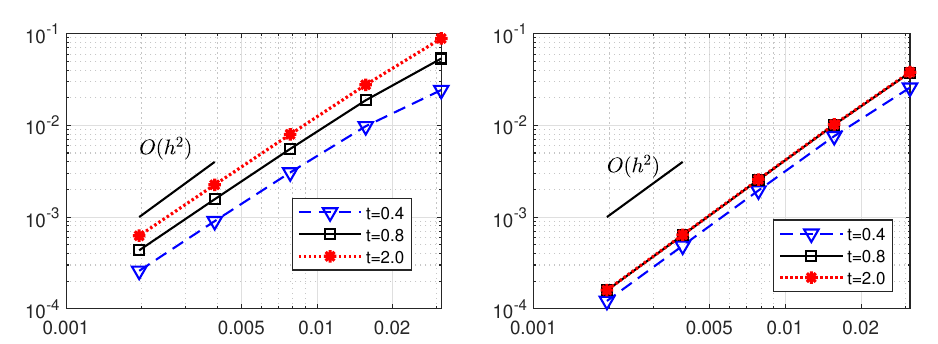}
	\caption{ [$\bkap=-2$] Numerical errors in the evolution of  an initial circle segment under Navier boundary conditions (left panel) and clamped boundary conditions (right panel),  where $\ttau=(\frac{2^{5}\,h}{5})^2 $. }
	\label{fig:openerror}
	\vspace{-0.3cm}
\end{figure}

\begin{figure}[!htp]
	\centering
	\includegraphics[width=0.9\textwidth]{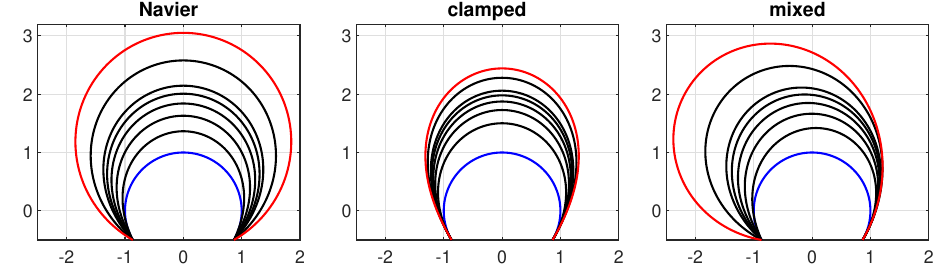}
	\includegraphics[width=0.9\textwidth]{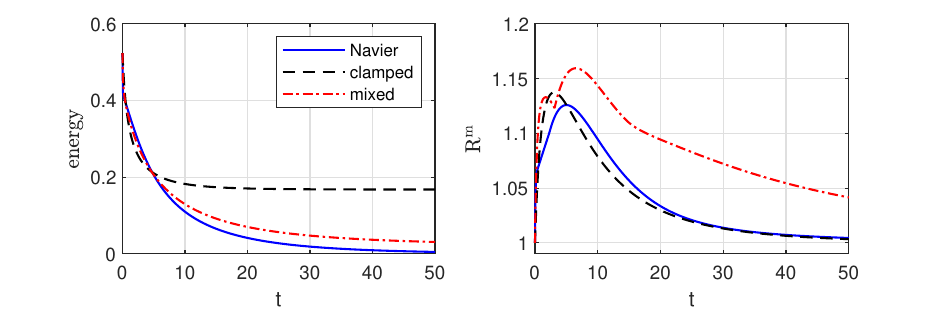}
	\caption{[$\bkap=-0.5$] Evolution of an initial circle segment towards the steady state (red line) under different boundary conditions. We plot $\Gamma^m$ at times $t=0,2,4, \cdots, 20, 50$. On the bottom are plots of the discrete energy and the mesh ratio ${\rm R}^m$. }
	\label{fig:bkap05bd}
\end{figure}

\vspace{0.3cm}
\noindent
{\bf  Example 2}:  We next perform convergence experiments for a circle segment of radius $r_0=1$ and central angle $\frac{2\pi}{3}$ under different boundary conditions, see Fig.\ \ref{fig:bkap05bd}. For an open curve $\Gamma_1$, we denote by  $\Omega_1$ the region enclosed by $\Gamma_1$ with the line segment connecting the two endpoints of $\Gamma_1$.  Similarly we introduce the region $\Omega_2$ for the open curve $\Gamma_2$.  Then the difference between $\Gamma_1$ and $\Gamma_2$ is measured by the area of the symmetric difference between $\Omega_1$ and $\Omega_2$ \cite{Zhao2021energy}
\begin{equation}
{\rm MD}(\Gamma_1,\Gamma_2) = |(\Omega_1\setminus\Omega_2)\cup(\Omega_2\setminus\Omega_1)| = |\Omega_1| + |\Omega_2|-2 |\Omega_1\cap\Omega_2|,\nn 
\end{equation}
where $|\Omega|$ represents the area of the region $\Omega$.  We then define the numerical errors
\begin{equation}
e_{h,\ttau}(t) = {\rm MD }(\Gamma_{h,\ttau}(t), \Gamma_{\frac{h}{2},\frac{\ttau}{4}}(t)), \nn 
\end{equation}
where $\Gamma^h_{h,\ttau}(t) = \vec X_{h,\ttau}(\Gamma^m, t)$ for $t\in [t_m, t_{m+1}]$ and $\vec X_{h,\ttau}(\cdot, t)$ is defined via \eqref{eq:linertime}.

  We fix $\bkap=-2$ and consider both the cases of Navier boundary and clamped boundary conditions. As expected,  a second-order convergence rate for the numerical solutions can be observed as well in Fig.~\ref{fig:openerror}. 

\begin{figure}[!htp]
	\centering
	\includegraphics[width=0.9\textwidth]{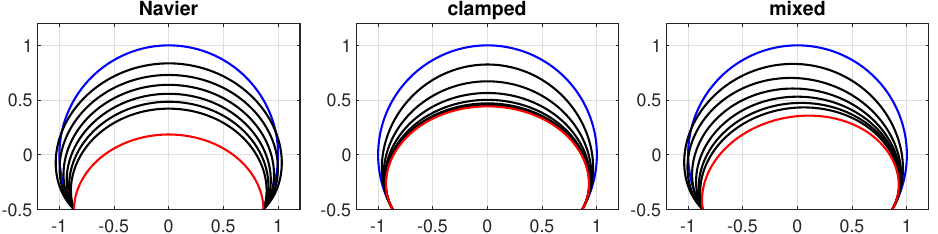}
	\includegraphics[width=0.9\textwidth]{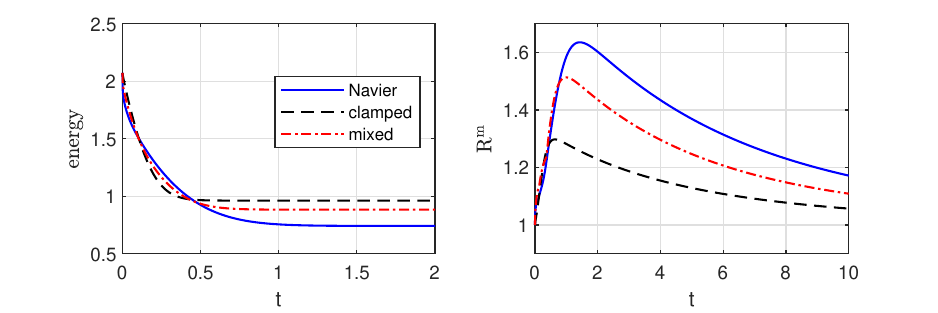}
	\caption{[$\bkap=-2$] Evolution of an initial circle segment towards the steady state (red line) under different boundary conditions. We plot $\Gamma^m$ at times $t=0,0.1,\cdots,0.6, 2$.  On the bottom are plots of the discrete energy and the mesh ratio ${\rm R}^m$. }
	\label{fig:bkap2bd}
\end{figure}

\begin{figure}[!htp]
	\centering
	\includegraphics[width=0.9\textwidth]{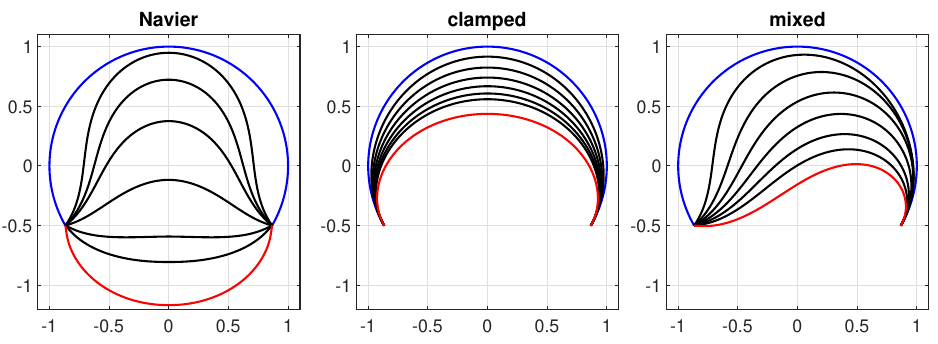}
	\includegraphics[width=0.9\textwidth]{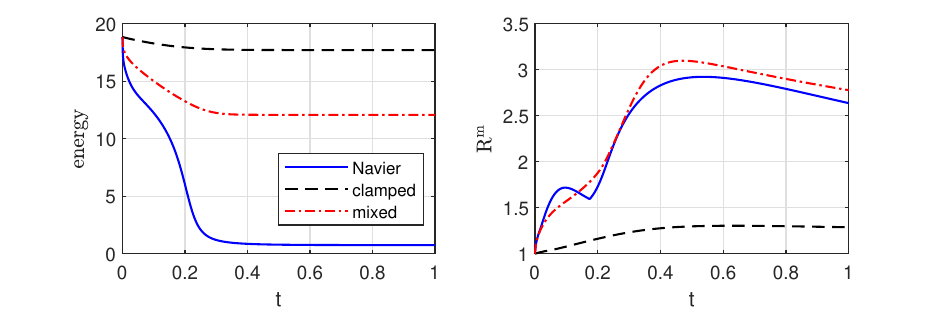}
	\caption{[$\bkap=2$] Evolution of an initial circle segment towards the steady state (red line) under different boundary conditions. We plot $\Gamma^m$ at times $t=0,0.05,\cdots,0.3, 1$.  On the bottom are plots of the discrete energy and the mesh ratio ${\rm R}^m$. }
	\label{fig:bkapp2bd}
\end{figure}

\vspace{0.3cm}
\noindent
{\bf Example 3}:  In this example, we consider the evolution of an open curve under different types of boundary conditions. We also monitor the mesh ratio
\begin{equation*}
{\rm R}^m = \frac{\max_{\sigma\in\mathscr{T}^m} \mH^{d-1}(\sigma)}{\min_{\sigma\in\mathscr{T}^m} \mH^{d-1}(\sigma)},\quad m=0,1,\ldots, M, 
\end{equation*} 
which implies that the vertices on the polygonal curves are equidistributed if ${\rm R}^m=1$. 

Here the initial setting is the same as that in {\bf  Example 2} for a circle segment. We fix $J=128$ and $\ttau= 10^{-3}$ and conduct an experiment with $\bkap=-0.5$, which will lead to the curve expanding.
 The numerical results are reported in Fig.~\ref{fig:bkap05bd},  where the effects of different boundary conditions are visualized.
Indeed, in the case of Navier boundary conditions, the curve evolves to a circle
segment of radius 2, whereas in the case of a clamped boundary, the tangent vectors at the boundary remain unchanged, which numerically confirms  \eqref{eq:clamptau}. Moreover, in the case of mixed boundary conditions the curve evolves towards an asymmetric shape.  We also observe the decay of the discrete energy in the three cases, which substantiates Theorem~\ref{thm:stability}.  Moreover, the mesh ratio function only grows to a very small value and then gradually decrease to 1. This implies the mesh quality is well preserved in our simulations.

We next perform an experiment with $\bkap=-2$ and the results are shown in Fig.~\ref{fig:bkap2bd}. Notice that due to the distance of the two boundary points, this time the curve is not able to attain the zero energy shape of a circle segment with radius $0.5$. Once again we observe the energy decay and good mesh properties for all three types of boundary conditions.  We further repeat the experiment with $\bkap=2$ and the results are shown in Fig.~\ref{fig:bkapp2bd}. We notice that the curve attempts to flip the sign of its curvature, which it is prevented from doing in the case of clamped boundary conditions, as this would cause a too large deformation.  

\subsection{3d results} In this subsection, we perform numerical tests for surfaces in $\bR^3$, both for closed surfaces and for surfaces with boundary with different types of boundary conditions. 

\begin{figure}[t]
	\centering
	\includegraphics[width=0.45\textwidth]{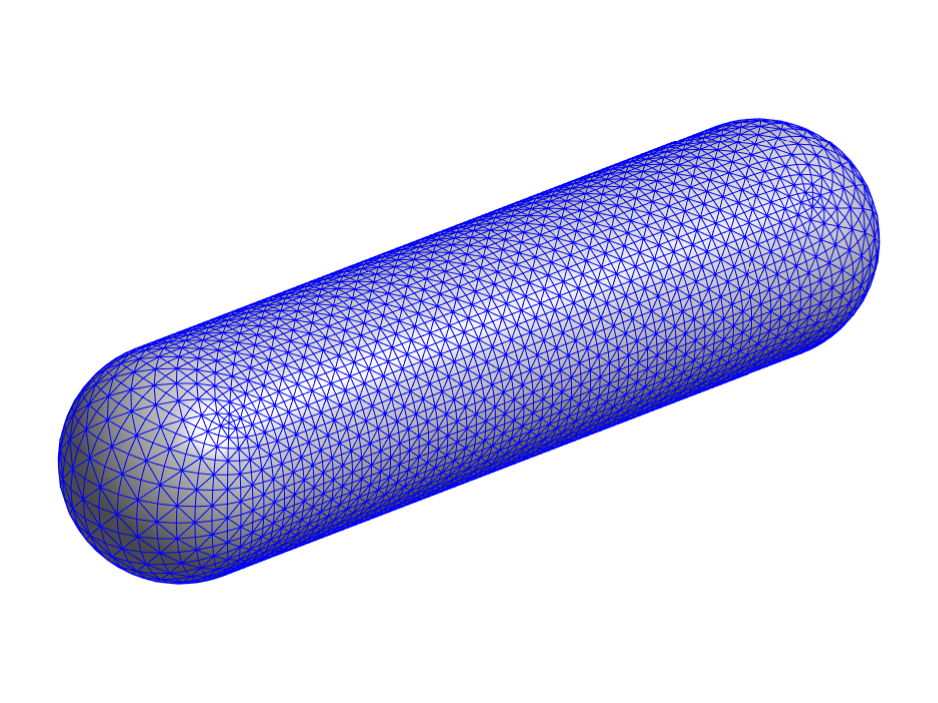}
	\includegraphics[width=0.45\textwidth]{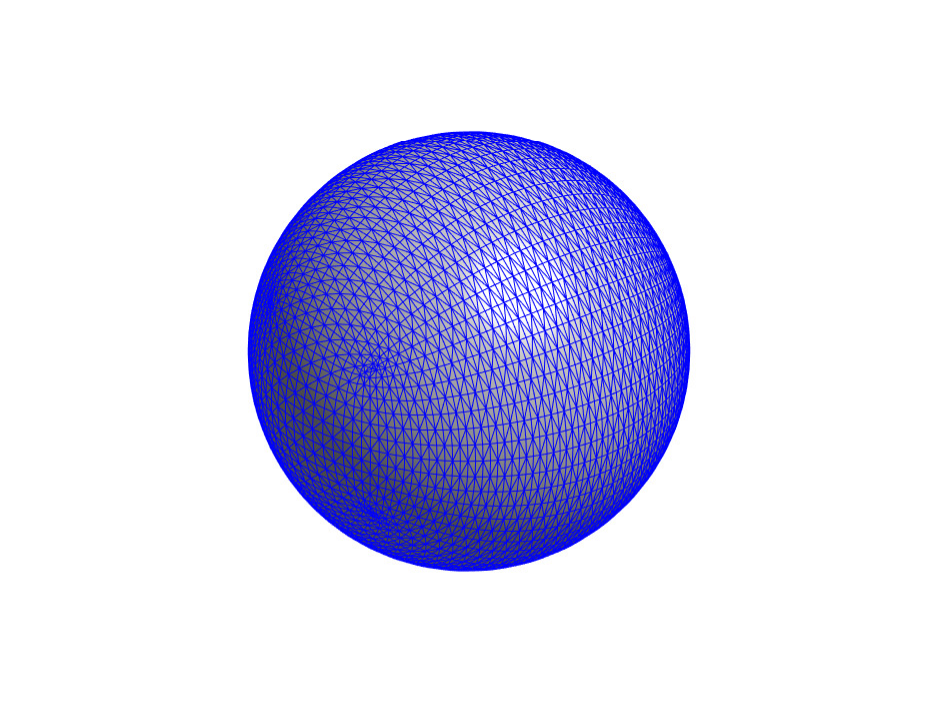}
	\includegraphics[width=0.7\textwidth]{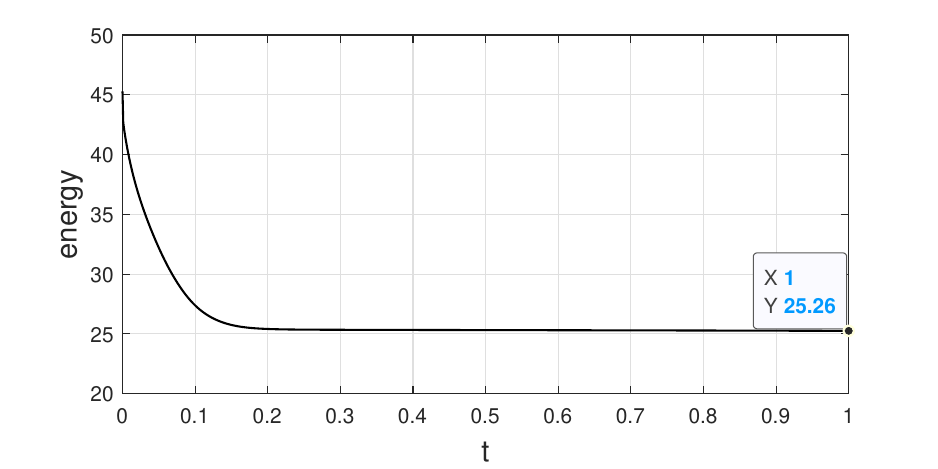}
	\caption{[$\bkap=0$] Evolution of an initial cigar shape of dimension $4\times 1\times 1$ towards the steady state, where we plot $\Gm$ at $t=0, 1$. On the bottom is a plot of the discrete energy.}
	\label{fig:411cigar}
\end{figure}

\vspace{0.3cm}
\noindent
{\bf Example 4}:  We begin with an experiment for a cigar shape of total dimension $4\times 1\times 1$. We fix $\ttau=10^{-3}$ and choose $\bkap=0$. Here $J=9216, K = 4610$.  The numerical results are shown in Fig.~\ref{fig:411cigar}. As expected, we observe that the round cylinder evolves towards a sphere as the steady state, where the final discrete energy is about 25.26, which approximates the value $8\pi=25.13$, the Willmore energy of a sphere. Here the energy decay property is observed and the mesh quality of the polyhedral surface is in generally well preserved during the evolution. 

\begin{figure}[t]
	\centering
	\includegraphics[width=0.44\textwidth]{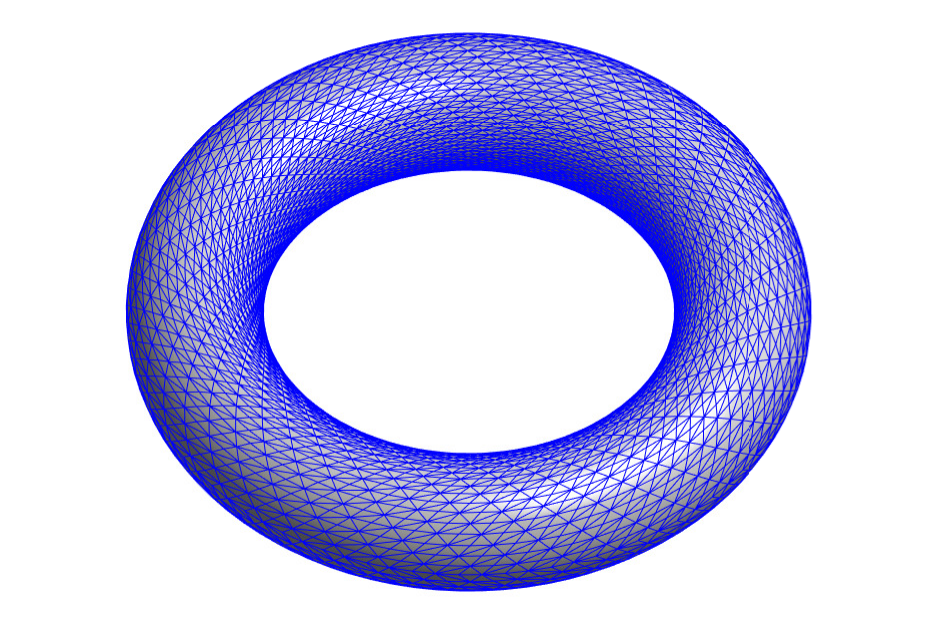}
	\includegraphics[width=0.44\textwidth]{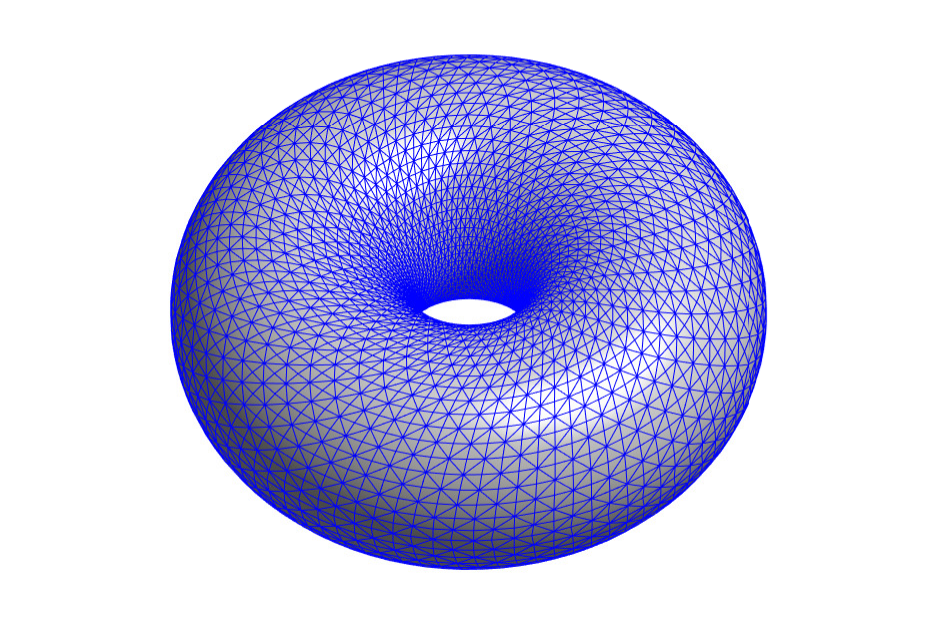}
	\includegraphics[width=0.7\textwidth]{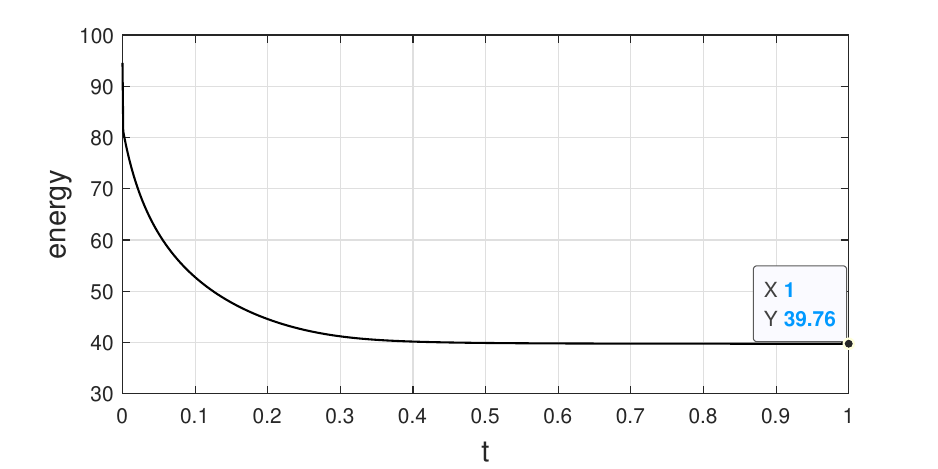}
	\caption{[$\bkap=0$] Evolution of an initial thin torus towards a Clifford torus. We plot $\Gm$ at $t=0, 1$.  On the bottom is a plot of the discrete energy.  }
	\label{fig:torus}
	\vspace{-0.5cm}
\end{figure}

We next consider an experiment for a thin torus with major radius $2$ and minor radius 0.5. We fix $\ttau=10^{-3}$, $\bkap=0$ and $J=11040, K = 5520$. The numerical results are reported in Fig.~\ref{fig:torus}, where we observe that the torus grows towards a Clifford torus as the steady state. In fact, the energy at the final time $T=1$ is $39.76$, which approximates the value $4\pi^2 = 39.48$, the Willmore energy of the Clifford torus. Recall that the Clifford torus has been proved to be the energy minimizer among surfaces of genus 1 in \cite{Marques14min}.

\begin{figure}[!htp]
	\centering
	\includegraphics[width=0.44\textwidth]{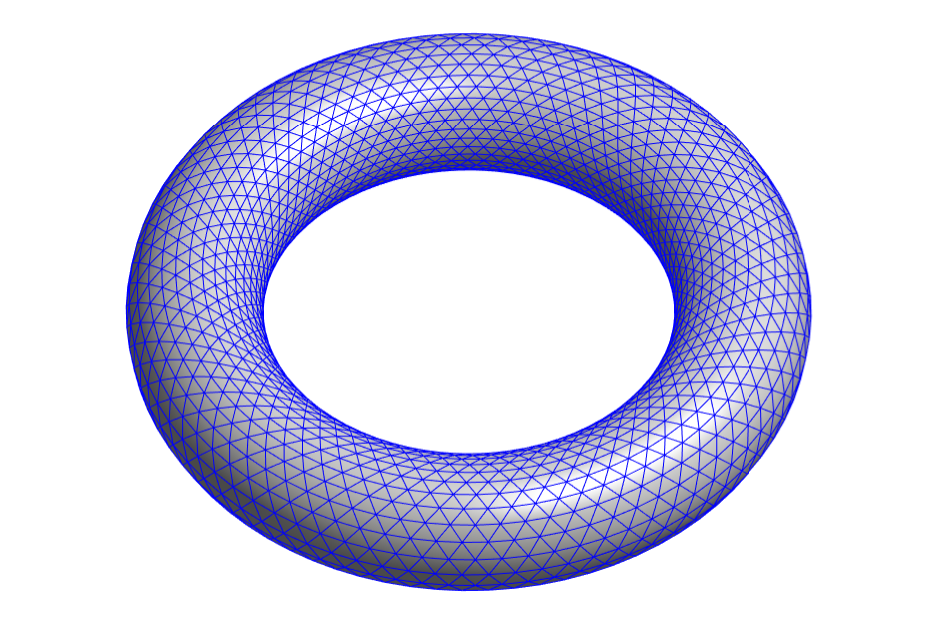}
	\includegraphics[width=0.44\textwidth]{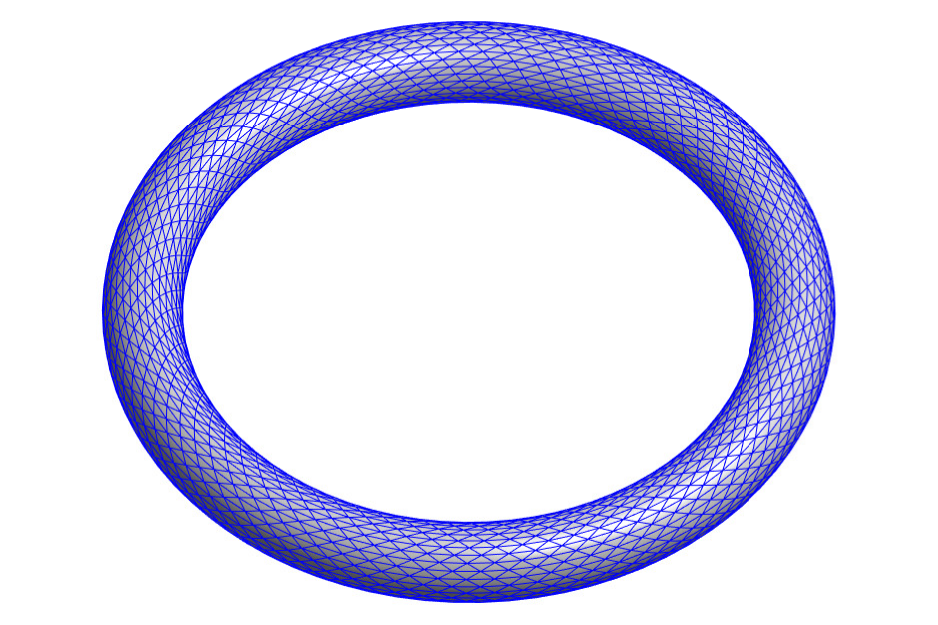}
	\includegraphics[width=0.7\textwidth]{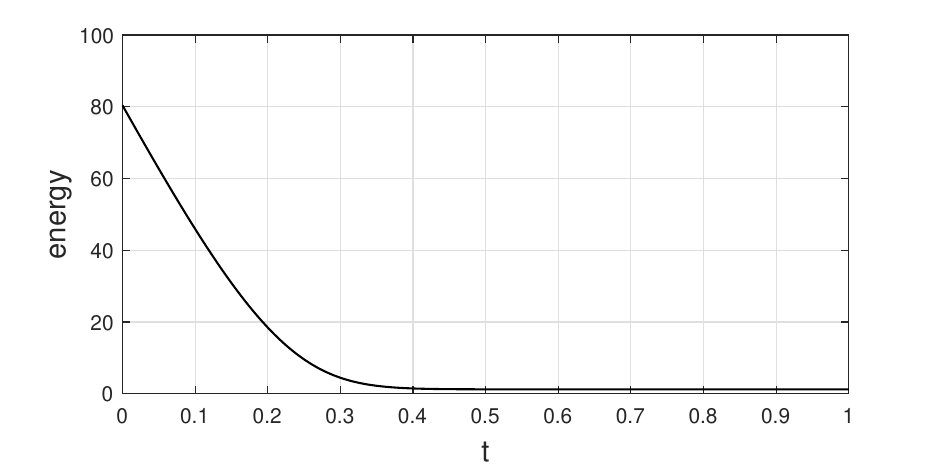}
	\caption{[$\bkap=-2$] Evolution for a torus, where we plot $\Gm$ at $t=0,1$. On the bottom is a plot of the discrete energy. }
	\label{fig:torusbkap-2}\vspace{-0.5cm}
\end{figure}

\begin{figure}[!htp]
	\centering
	\includegraphics[width=0.45\textwidth]{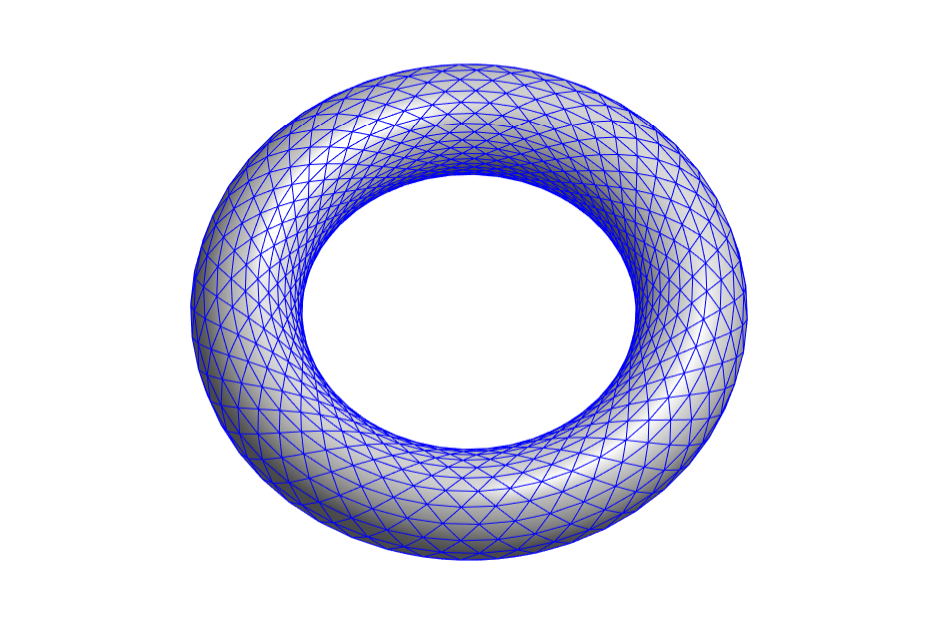}
	\includegraphics[width=0.45\textwidth]{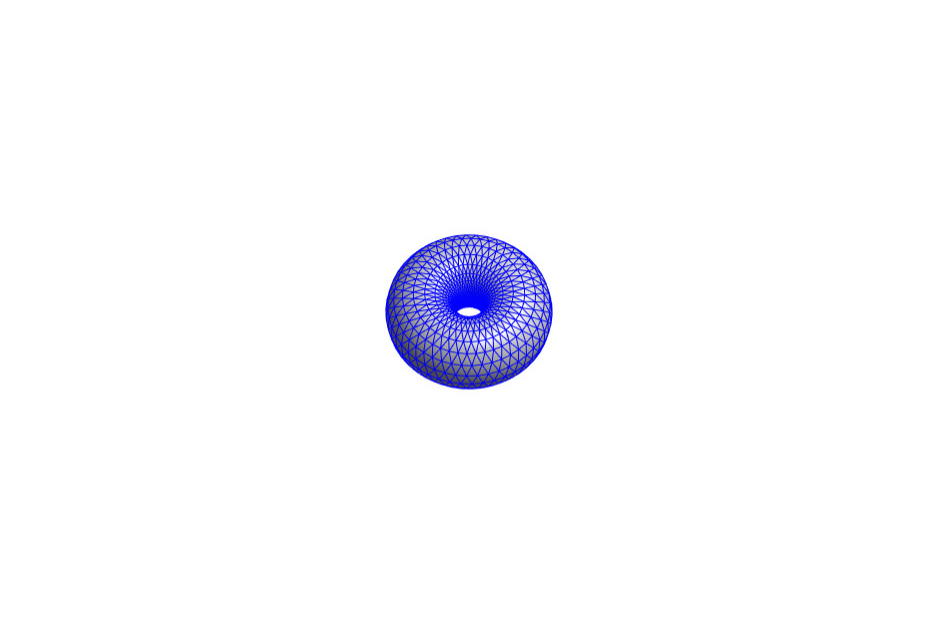}
	\includegraphics[width=0.7\textwidth]{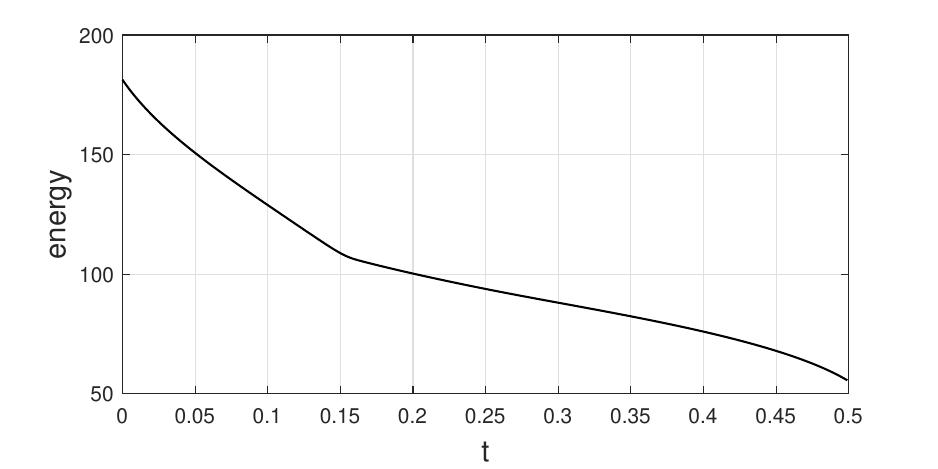}
	\caption{[$\bkap=1$] Evolution for a torus, where we  plot $\Gm$ at $t=0, 0.5$.  On the bottom is a plot of the discrete energy.  }
	\label{fig:torusbkap1}
\end{figure}

We also conduct experiments for tori with $\bkap=-2$ and $\bkap=1$, respectively. The numerical results are shown in Fig.~\ref{fig:torusbkap-2} and Fig.~\ref{fig:torusbkap1}. Here we observe that the major radius of the torus either grows or decreases in order to decrease the energy.

\begin{figure}[!htp]
	\centering
	\includegraphics[width=0.9\textwidth]{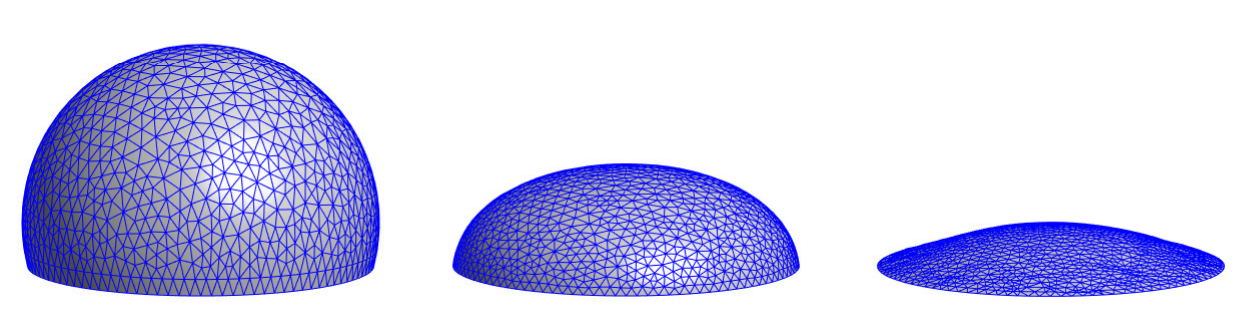}
	\includegraphics[width=0.9\textwidth]{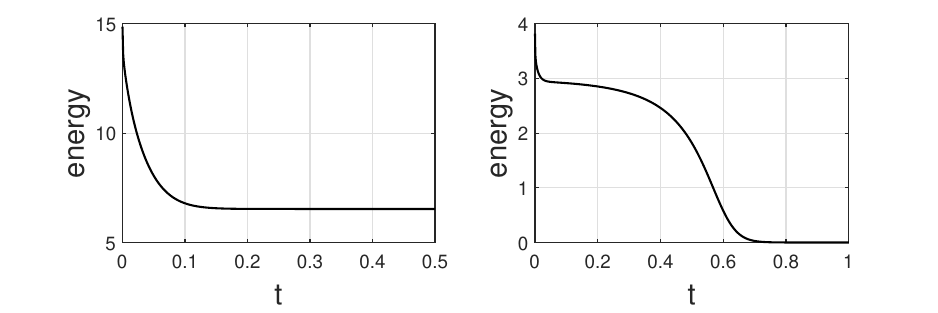}
	\caption{Evolution of an initial sphere cap under Navier boundary conditions. Left panel: $t=0$; middle panel:$\bkap=-4$ at $t=0.5$; right panel: $\bkap=-1$ at $t=1$. On the bottom is a plot of the discrete energy.}
	\label{fig:navierbd}
\end{figure}

\begin{figure}[!htp]
	\centering
	\includegraphics[width=0.9\textwidth]{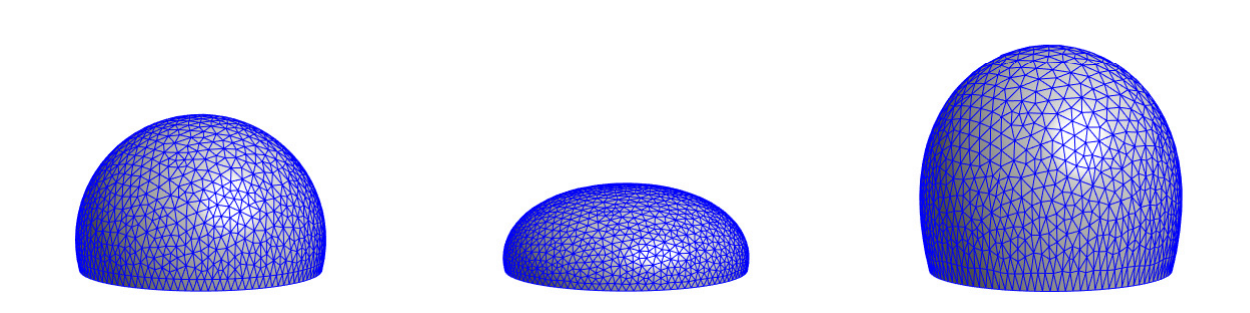}
	\includegraphics[width=0.9\textwidth]{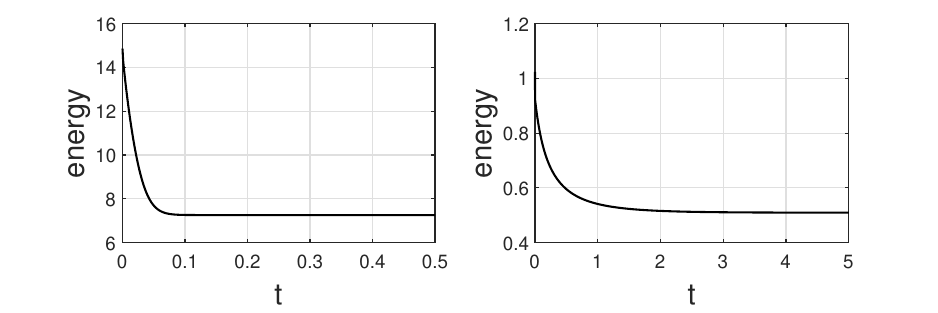}
	\caption{Evolution of an initial sphere cap under clamped boundary conditions. Left panel: $t=0$; middle panel: $\bkap=-4$ at $t=0.1$; right panel: $\bkap=-1.5$ at $t=5$.  On the bottom is a plot of the discrete energy. }
	\label{fig:clampbd}
\end{figure}

\vspace{0.3cm}
\noindent
{\bf Example 5}: We first consider the evolution of a sphere cap of radius 1 under Navier boundary conditions. We fix $\ttau=10^{-3}$ and use $J=2274, K=1188$.  The numerical results are demonstrated in Fig.~\ref{fig:navierbd} with $\bkap=-4$ and $\bkap=-1$.  We also repeat the above experiment under clamped boundary conditions with $\bkap=-4$ and $\bkap=-1.5$, and the results are visualized in Fig.~\ref{fig:clampbd}. Here we can clearly observe the effects of the spontaneous curvature and the boundary conditions on the evolution of the open surfaces. 

\begin{figure}[!htp]
	\centering
	\includegraphics[width=0.85\textwidth]{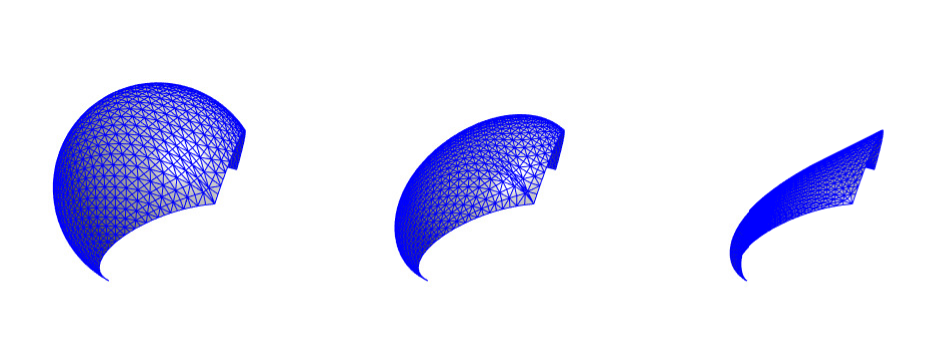}
	\includegraphics[width=0.85\textwidth]{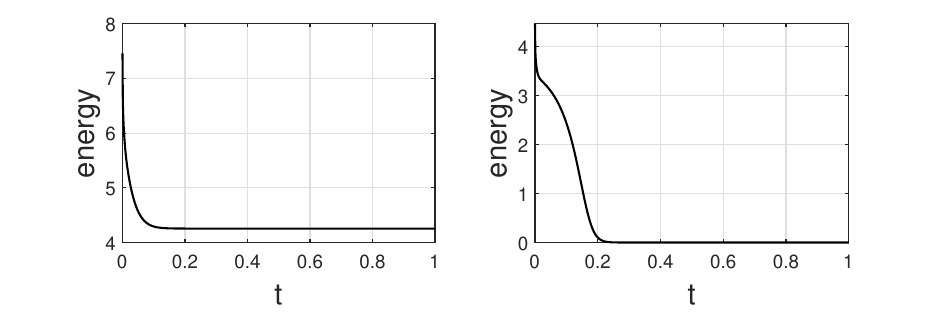}
		\caption{Evolution of an initial sphere cap from a standard quadruple bubble with Navier boundary conditions. Left panel: $t=0$; middle panel: $\bkap=-4$ at $t=1$; right panel: $\bkap=-1$ at $t=1$. On the bottom is a plot of the discrete energy.}
		\label{fig:nvqbd}
\end{figure}

\begin{figure}[!htp]
	\centering
	\includegraphics[width=0.85\textwidth]{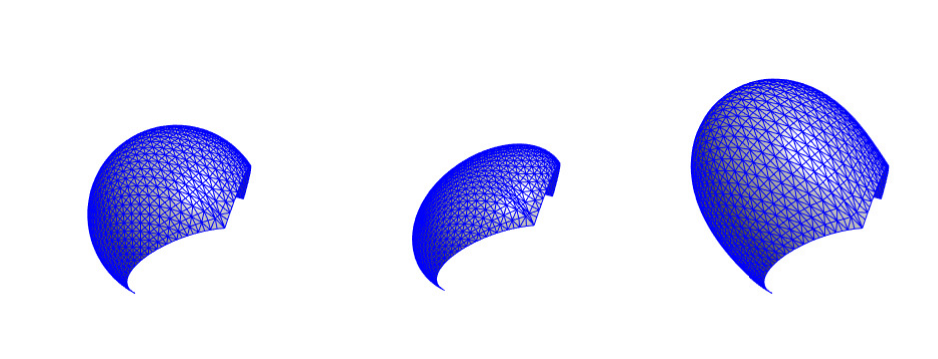}
	\includegraphics[width=0.85\textwidth]{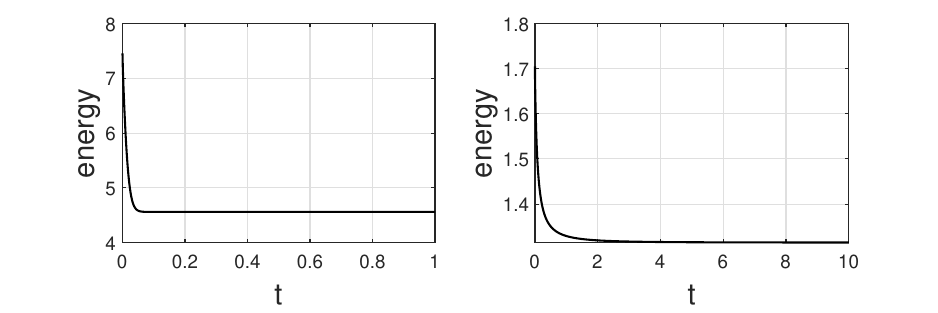}
\caption{Evolution of an initial spherical cap from a standard quadruple bubble with clamped boundary conditions. Left panel: $t=0$; middle panel: $\bkap=-4$ at $t=1$; right panel: $\bkap=-1.5$ at $t=10$. On the bottom is a plot of the discrete energy.}
\label{fig:clamqbd}
\end{figure}

\vspace{0.3cm}
\noindent
{\bf Example 6}:  We next include experiments where the boundary of the surface does not lie in a hyperplane. In fact, we take as initial data a sphere cap from a standard quadruple bubble, see \cite[Fig. 19]{cluster3d}. We fix $\ttau=10^{-3}$ and use $J=3072, K = 1601$.  Two experiments for Navier boundary conditions and clamped boundary conditions are shown in Fig.~\ref{fig:nvqbd} and Fig.~\ref{fig:clamqbd}, respectively.

\begin{figure}[t]
	\centering
	\includegraphics[width=0.49\textwidth]{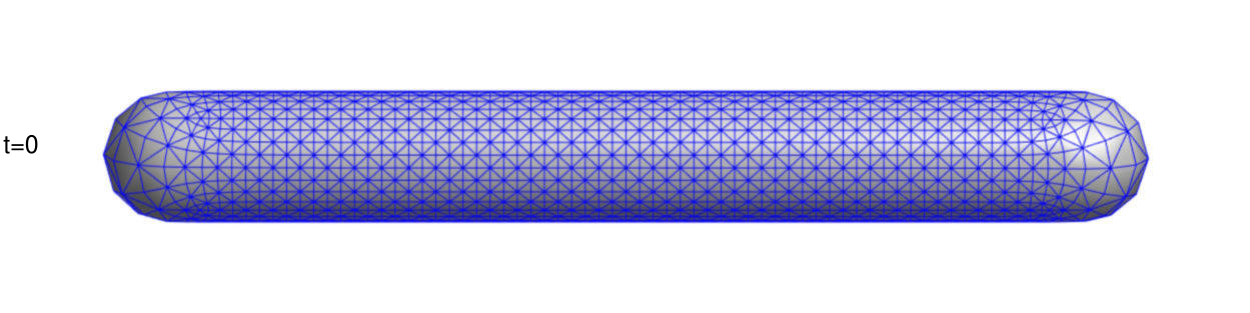}
	\includegraphics[width=0.49\textwidth]{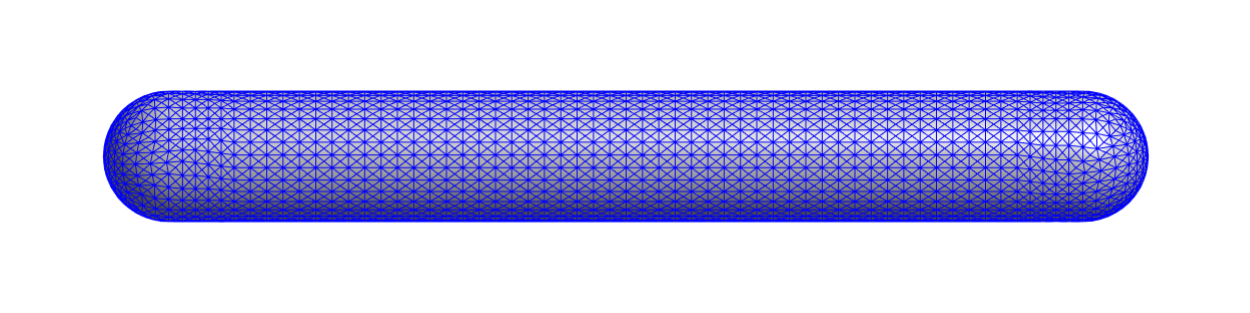}
	\includegraphics[width=0.49\textwidth]{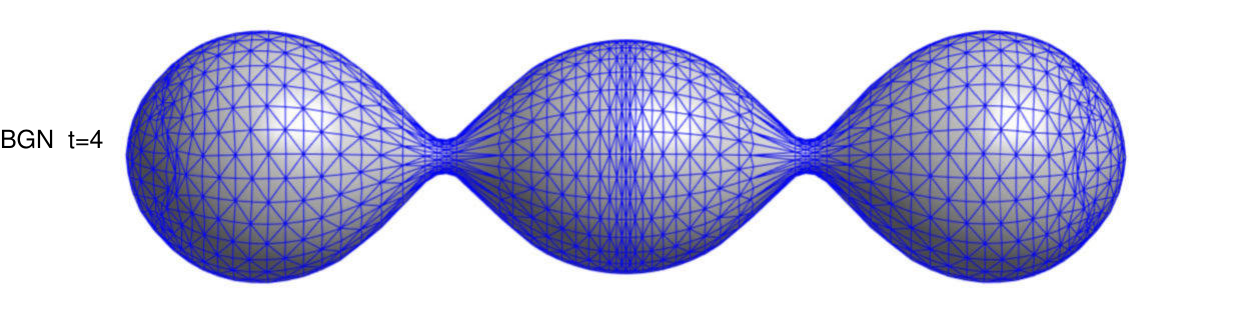}
	\includegraphics[width=0.49\textwidth]{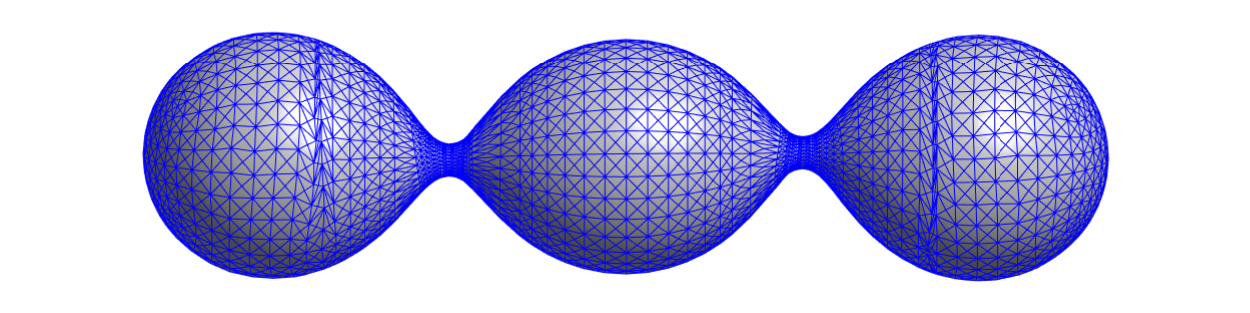}
	\includegraphics[width=0.49\textwidth]{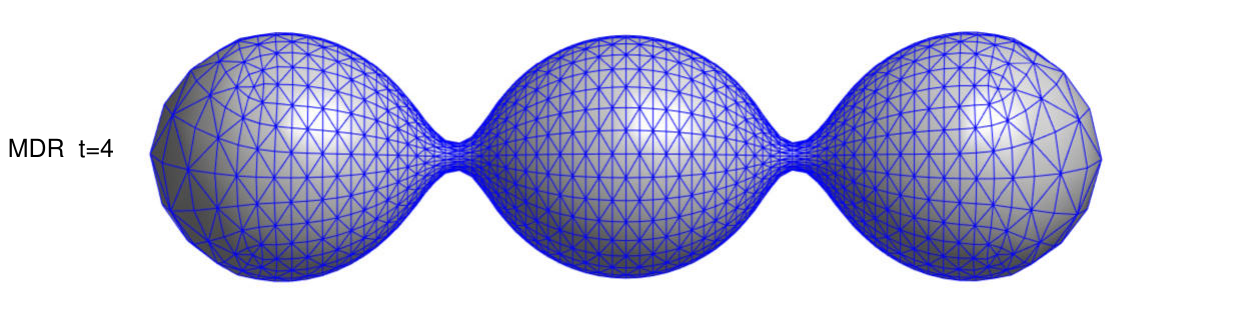}
	\includegraphics[width=0.49\textwidth]{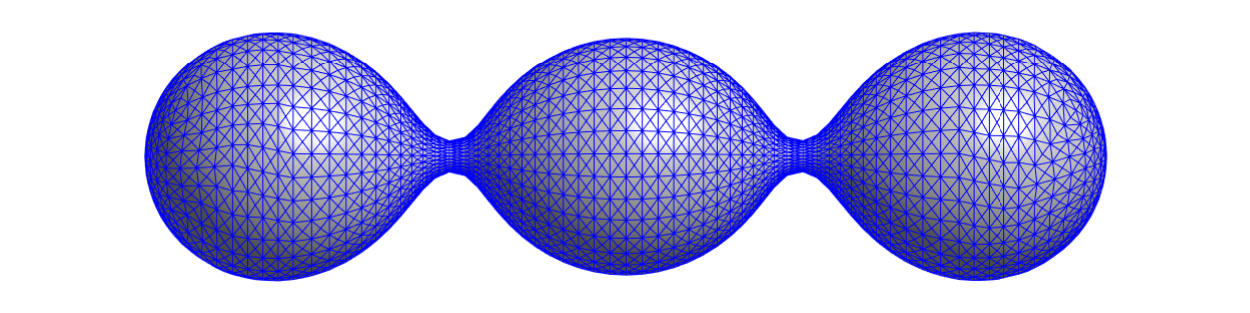}
	\includegraphics[width=0.49\textwidth]{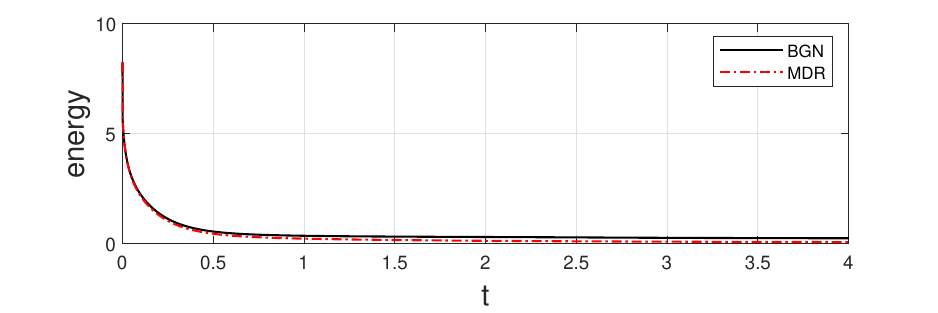}
	\includegraphics[width=0.49\textwidth]{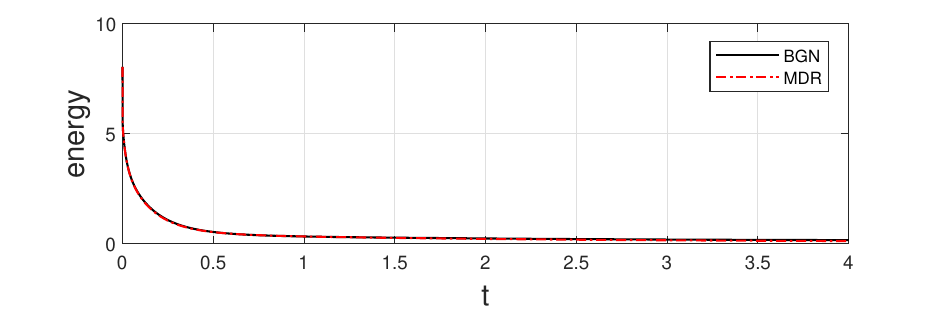}
	\caption{ [$\bkap=-2$] Evolution of an initial cigar shape of dimension of $8\times 1\times 1$ using the `BGN' or `MDR' approaches with a fixed time step size $\ttau= 10^{-3}$. Left panel: $J=4352, K =2178$; right panel: $J=8704, K = 4354$.}
	\label{fig:811cigar}
\end{figure}

\vspace{0.3cm}
\noindent
{\bf Example 7}: We end this section by considering an experiment for an initial cigar shape of total dimension $8\times 1\times 1$. We choose $\bkap=-2$ and use $\ttau=10^{-3}$  with two different computational surface meshes. The numerical results are shown in Fig.~\ref{fig:811cigar}, where we observe an evolution towards three connected spheres. Here we have employed tangential velocities based on the `BGN' approach in \eqref{eqn:fdVt} or the `MDR' approach that was considered in \eqref{eqn:fdVtmdr}. In general, both approaches can produce very similar results, which demonstrates the flexibility of our method in choosing different tangential velocities for the flow. In particular, we observe that the energy difference between the two approaches get smaller when using a finer computational mesh.

\section{Conclusions}\label{sec:con}

We proposed a new parametric finite element method for the Willmore flow of hypersurfaces in a unified framework. The method is linear and employs a splitting of the normal and tangential velocities of the flow. The normal velocity is approximated via an evolution equation for the curvature, and follows the arbitrary Lagrangian-Eulerian (ALE) approach. This enables an unconditionally stable scheme with respect to the discrete energy. We also incorporated the `BGN' tangential velocity into the flow through a curvature identity. This helps to preserve the mesh quality for the evolving discrete surface. We rigorously proved that the proposed method admits a unique solution and is energy stable. A variety of numerical examples were presented to verify the numerical convergence and demonstrate the favorable properties of the method.

The major advantages of our method lie in its unconditional energy stability, the linear structure, and flexibility in handling both closed and open hypersurfaces with complex boundary conditions. By decoupling the normal and tangential velocities, our method not only simplifies the numerical treatment but also allows for the independent selection of tangential velocities to improve mesh quality without sacrificing the stability. This offers a unified and practical computational parametric framework for numerically simulating Willmore flow.  

Extending the presented techniques to high-order parametric finite element spaces is not difficult. However, designing high-order temporal discretizations that preserve the energy decay property is nontrivial and requires further investigation.

\section*{Acknowledgement}
This work was partially supported by the National Natural
Science Foundation of China No. 12401572 (Q.Z).

\appendix 
\section{Differential calculus}\label{sec:appA}
\numberwithin{equation}{section}

We recall the integration by parts formula on the hypersurface $\Gt$ as
\begin{equation}
	\int_{\Gamma}\nabs\phi\dH^{d-1} =- \int_{\Gamma}\phi\,\varkappa\,\vec\nu\dH^{d-1}+\int_{\partial\Gamma}\phi\,\vec\mu\dH^{d-2},\label{eq:interpart}
\end{equation}
where $\varkappa$ is the mean curvature of $\Gamma$, $\vec\nu$ is the unit normal to $\Gamma$, and $\vec\mu$ is the outer unit conormal to $\partial\Gamma$.

Given $f\in C^1(\mathcal{G}_T)$, the Reynolds transport theorem reads as (see \cite[Theorem 32]{Barrett20})
\begin{align}
	\ddt\int_{\Gt} f\dH^{d-1} &= \int_{\Gt}(\partial_t^\circ f + f\nabs\cdot\mathscr{\vv V}\dH^{d-1} \nn\\
	&=\int_{\Gt}\left(\partial_t^\square f - f\,\mathscr{V}\,\varkappa\right)\dH^{d-1} + \int_{\partial\Gt} f\,\mathscr{\vv V}\cdot\vec\mu\dH^{d-2},\label{eq:transport}
\end{align}
where the second equality results from integration by parts.

\section{Proof of \eqref{eq:nstab2}}\label{sec:nstab2}
It holds for any matrices $\mat{A}=(a_{ij})\in\bR^{d\times d}$ and $\mat{B}=(b_{ij})\in\bR^{d\times d}$ that
\[\mat{A}:(\mat{A}-\mat{B}) = \sum_{i,j=1}^d a_{ij}(a_{ij}-b_{ij})\geq \frac{1}{2}\sum_{i,j=1}^d (a_{ij}^2-b_{ij}^2)=\frac{1}{2}(|\mat{A}|^2-|\mat{B}|^2).\]
Then we can compute 
\begin{align}
&\ipd{\nabs\vec X^{m+1}: \nabs[\vec X^{m+1}-\vec\id], [\varkappa^{m+1}-\bkap]^2}_{\Gm}\nn\\&\qquad
\geq \frac{1}{2}\ipd{|\nabla_s\vec X^{m+1}|^2 - |\nabla_s\vec\id|^2,~[\varkappa^{m+1}-\bkap]^2}_{\Gamma^m}\nn\\ & \qquad
=\frac{1}{2}\ipd{|\nabla_s\vec X^{m+1}|^2 - (d-1),~[\varkappa^{m+1}-\bkap]^2}_{\Gamma^m},\label{eq:astab1}
\end{align}
where we used the fact that $|\nabla_s\vec\id|^2=d-1$. 
     
In the case $d=2$, $\Gamma^m$ is a polyhedral curve. Then we can use the inequality $\frac{1}{2}(a^2-1)\geq (a-1)$ to obtain that 
\begin{align}
&\frac{1}{2}\ipd{|\vec X_s^{m+1}|^2 - 1,~[\varkappa^{m+1}-\bkap]^2}_{\Gamma^m}\geq \ipd{|\vec X_s^{m+1}|-1,~[\varkappa^{m+1}-\bkap]^2}_{\Gamma^m}\nn\\
&\qquad\qquad  = \ipd{1,~[\varkappa_{\Gamma^{m+1}}^{m+1}-\bkap]^2}_{\Gamma^{m+1}}-\ipd{1,~[\varkappa^{m+1}-\bkap]^2}_{\Gamma^m}.\label{eq:astab2}
\end{align}

In the case $d=3$, we can compute the last term in \eqref{eq:astab1} as
\begin{align}
&\frac{1}{2}\ipd{|\nabla_s\vec X^{m+1}|^2 - (d-1),~[\varkappa^{m+1}-\bkap]^2}_{\Gamma^m}\nn\\
&=\frac{1}{2}\ipd{|\nabla_s\vec X^{m+1}|^2,~[\varkappa^{m+1}-\bkap]^2}_{\Gamma^m} -\ipd{1, [\varkappa^{m+1}-\bkap]^2}_{\Gamma^m}.\label{eq:astab3}
\end{align}
On each $\sigma\in\mT^m$, we now recall the inequality  (see \cite[(26)]{Barrett20})
\begin{equation}
\sqrt{g} \leq \frac{1}{2}\left(|\nabla_s\vec X^{m+1}\vec t_1|^2 +|\nabla_s\vec X^{m+1}\vec t_2|^2\right)=\frac{1}{2}|\nabla_s\vec X^{m+1}|^2,
\end{equation}
where $\{\vec t_1, \vec t_2\}$ is an orthonormal basis for the tangent plane of $\sigma\in\mT^m$ and $g = {\rm det}([\nabla_s\vec X^{m+1}]^T\nabla_s\vec X^{m+1})$, which measures the square of the local area change of the mapping $\vec X^{m+1}: \Gamma^m\mapsto\Gamma^{m+1}$.  
Therefore we can recast the first term in \eqref{eq:astab3} as
\begin{align}
&\frac{1}{2}\ipd{|\nabla_s\vec X^{m+1}|^2,~[\varkappa^{m+1}-\bkap]^2}_{\Gamma^m}=\frac{1}{2}\sum_{j=1}^J\int_{\sigma_j^m}|\nabla_s\vec X^{m+1}|^2 [\varkappa^{m+1}-\bkap]^2\,\dH^{d-1}\nn\\
&\qquad\qquad\geq \sum_{j=1}^J\int_{\sigma_j^m}\sqrt{g}\,[\varkappa^{m+1}-\bkap]^2\,\dH^{d-1}\nn\\
&\qquad\qquad = \sum_{j=1}^J\int_{\sigma_j^{m+1}} [\varkappa^{m+1}_{\Gamma^{m+1}}-\bkap]^2\,\dH^{d-1} = \norm{\varkappa^{m+1}_{\Gamma^{m+1}}-\bkap}_{m+1}^2.\label{eq:astab4}
\end{align}
Combining \eqref{eq:astab2} and \eqref{eq:astab4} with \eqref{eq:astab1} yields \eqref{eq:nstab2}. 
    
\bibliographystyle{siam}
\normalem
\bibliography{bib}

\end{document}